\documentclass[11pt]{article}
\usepackage[utf8]{inputenc}
\usepackage[T1]{fontenc}
\usepackage{amsmath,amssymb,amsthm,mathtools,xcolor,hyperref,listings,cleveref,tikz,centernot}
\usepackage{amsfonts}
\usepackage{mathrsfs}
\usetikzlibrary{arrows,decorations.pathreplacing,shapes,calc}
\usepackage{amscd}
\usepackage{enumitem}
\usepackage{geometry}
\usepackage{array}
\usepackage[authoryear]{natbib}
\bibpunct{(}{)}{,}{a}{}{,}
\usepackage{tcolorbox}
\tcbset{colback=gray!10, colframe=white, boxrule=0pt, arc=0.1mm, left=3mm, right=3mm, top=3mm, bottom=3mm}

\usepackage{colortbl}

\usepackage{titling}
\pretitle{\raggedright\LARGE\bfseries}
\posttitle{\par\vspace{1em}}

\preauthor{\raggedright\large}
\postauthor{\par}

\hypersetup{
      colorlinks=true,
      linkcolor=magenta,
      citecolor=teal,
      urlcolor=magenta,
      pdftitle={Non-Hausdorff Incidence Completions of Finite Coverings: Monodromy, Transport, and Defect Posets},
      pdfauthor={Abhiram Sripat},
      pdfkeywords={non-Hausdorff topology, finite covering spaces, exceptional fibres, peripheral monodromy, incidence, finite topological spaces, relation-valued transport, defect posets, Grothendieck construction, Artin-Wraith glueing}
}
\newcommand{\Sym}{\mathfrak{S}}

\theoremstyle{plain}
\newtheorem{theorem}{Theorem}[section]
\newtheorem{proposition}[theorem]{Proposition}
\newtheorem{corollary}[theorem]{Corollary}
\newtheorem{lemma}[theorem]{Lemma}

\theoremstyle{definition}
\newtheorem{definition}[theorem]{Definition}
\newtheorem{example}[theorem]{Example}
\newtheorem{remark}[theorem]{Remark}

\title{Non-Hausdorff Incidence Completions of Finite Coverings: Monodromy, Transport, and Defect Posets}
\author{
 Abhiram Sripat\\
 \textit{Florence Quantum Labs, California, USA}\\
 \href{mailto:abhiram@florencequantumlabs.com}{\textit{abhiram@florencequantumlabs.com}}
}
\date{}

\begin{document}
\pagestyle{empty}
\maketitle
\thispagestyle{empty}

\begin{tcolorbox}
\noindent\textbf{Abstract.}
Let $X$ be a connected Hausdorff surface, let $\Sigma\subset X$ be finite, and let
$\pi^\circ:Y^\circ\to X\setminus\Sigma$ be a finite-sheeted covering.  For each
marked point, peripheral monodromy determines a finite set of peripheral component
germs $C_i$ with local degree labels $\nu_i$.  We prescribe a topological exceptional
fibre $F_i$ and nonempty closed endpoint sets $A_{i,c}\subseteq F_i$, allowing one
germ to have several exceptional endpoints and several germs to share an endpoint.

The associated tail-saturated topology recovers the prescribed incidence intrinsically:
the endpoint set of the $c$-tail is exactly $A_{i,c}$.  A fixed regular lift accumulating
at a defect through the germ set $\Lambda$ extends through that point exactly when its
exceptional value lies in $\bigcap_{c\in\Lambda}A_{i,c}$.  For finite discrete fibres,
the bipartite incidence graph $G_{\mathrm{inc}}$ controls the local fundamental group
and low-dimensional homology, with
$H_1\cong H_1(G_{\mathrm{inc}})$ and
$H_2\cong\bigoplus_c\widetilde H_0(A_c)$.  Power and multi-arm probes recover the
peripheral degrees and the discrete incidence relation by M\"obius inversion.

For finite exceptional fibres, the completion datum is equivalently a weighted
specialization poset $(P,\mu)$ with an antitone incidence map
$I:P^{\mathrm{op}}\to\mathcal P(C)\setminus\{\varnothing\}$.  These data give a
structured monodromy classification, a singular automorphism kernel, and a structured
quotient correspondence.  With fixed local crossing data, incidence also determines
relation- and multiplicity-valued transport through isolated defects.  Finally, varying
the finite exceptional topology and incidence gives a Grothendieck construction over
the poset of finite topologies.  A full-incidence closure retraction shows that ordinary
total order-complex homotopy does not detect the incidence fibres.  In the one-end
sector, the classical Björner--Welker complex connects the top-homology representation
to partition-lattice/free-Lie representation theory; incidence-dependent information is
retained by the vertical fibres and their survival under topology degeneration.
\end{tcolorbox}

\noindent\textbf{Keywords:}
non-Hausdorff topology; finite covering spaces; exceptional fibres; peripheral monodromy; incidence; finite topological spaces; relation-valued transport; defect posets; Grothendieck construction; Artin--Wraith glueing.

\vspace{1em}
\sloppy
\pagestyle{plain}
\thispagestyle{empty}
\tableofcontents
%----------------------------------------------------------------------

\section{Introduction}
\label{sec:introduction}

Let $X$ be a connected Hausdorff surface without boundary, let
$\Sigma=\{z_1,\ldots,z_n\}\subset X$ be finite, and let
$\pi^\circ:Y^\circ\to X^\circ:=X\setminus\Sigma$ be a finite-sheeted covering.
For a sufficiently small punctured disc about $z_i$, its inverse image has a finite
set $C_i$ of connected components.  Standard covering theory identifies these
components with the cycles of peripheral monodromy, whose lengths define
$\nu_i:C_i\to\mathbb N_{>0}$ \citep{Hatcher2002,Spanier1966}.

We prescribe a topological exceptional fibre $F_i$ together with nonempty closed
endpoint sets $A_{i,c}\subseteq F_i$ for $c\in C_i$, with
$F_i=\bigcup_{c\in C_i}A_{i,c}$.  Equivalently, the datum is an incidence relation
$R_i\subseteq C_i\times F_i$ defined by $(c,a)\in R_i$ if and only if
$a\in A_{i,c}$.  A peripheral component germ may have several exceptional endpoints,
and several germs may share one endpoint.  Point replacement, Artin--Wraith glueing, and
Fox-type singular-covering completion provide the relevant classical background
\citep{Fox1957,Montesinos2005,Niefield1982,Niefield2012,deSouza2023}.  At the level of topological glueing, the construction below lies within this general
setting.  The additional structure is the incidence between peripheral component germs
of a genuine punctured covering and exceptional points, together with the lifting,
homological, reconstruction, symmetry, transport, and degeneration invariants derived
from that incidence.  We do not impose Fox's completeness or unique-completeness
conditions; in particular, a single peripheral component germ may have several
prescribed exceptional endpoints.  The term \emph{incidence completion} below refers to
this adjoining construction and does not assert completeness in Fox's sense.

The tail-saturated completion topology realizes the prescribed endpoint sets exactly,
\[
\operatorname{Ends}(c)=A_{i,c},
\]
and a fixed regular lift extends through a defect point $z$ exactly when its
exceptional value lies in
\[
\bigcap_{c\in\Lambda(z)}A_{i,c},
\]
where $\Lambda(z)$ is the set of accumulated peripheral ends.  For finite
discrete fibres the incidence graph $G_{\mathrm{inc}}=(C\sqcup F,E)$ controls
the local topology:
\[
H_1(\mathcal N;\mathbb Z)\cong H_1(G_{\mathrm{inc}};\mathbb Z),
\qquad
H_2(\mathcal N;\mathbb Z)\cong
\bigoplus_{c\in C}\widetilde H_0(A_c;\mathbb Z),
\]
and each connected component has free fundamental group of rank equal to the
cycle rank of the corresponding component of $G_{\mathrm{inc}}$.  Based power
probes recover $\nu(c)$ and endpoint multiplicity, while multi-arm probes
recover all intersection numbers
$\lambda(J)=|\bigcap_{c\in J}A_c|$; Boolean M\"obius inversion reconstructs
the discrete incidence graph.

For a finite exceptional fibre, classical finite-space theory
\citep{Alexandroff1937,Stong1966,McCord1966,Barmak2011} reduces the topology to
its Kolmogorov specialization poset $P$ and multiplicity function $\mu$; the
closed endpoint sets are equivalently an antitone map
\[
I:P^{\mathrm{op}}\longrightarrow\mathcal P(C)\setminus\{\varnothing\}.
\]
This yields structured monodromy data, a categorical classification of the
completions considered here, a singular automorphism kernel, and a structured
quotient correspondence.  Incidence also induces noninvertible transport
through defects, valued in finite relations and, for discrete fibres, matrices
over $\mathbb N$.

Finally, varying the finite exceptional topology and incidence produces the finite
defect poset $\mathfrak D_{k,r}$, the Grothendieck construction of the incidence fibres
over the finite-topology poset.  A compatible full-incidence maximum in every fibre
induces a closure retraction, so ordinary total order-complex homotopy is blind to the
incidence direction.  The one-end sector is the classical finite-topology/preorder
complex and is related to partition-lattice/free-Lie representation theory; for $r>1$,
incidence-dependent information is retained by the vertical incidence fibres and their
survival under topology degeneration.

\section{Multiple-origin lines}
\label{sec:multiple-origin}

The line with multiple origins supplies the minimal example in which distinct
points complete the same regular local germ without defining a covering projection.

Throughout this section let $k\geq 2$.

\subsection{The line with \texorpdfstring{$k$}{k} origins}
\label{subsec:Lk-definition}

Let
\[
\bigsqcup_{i=1}^k \mathbb R_i
\]
be the disjoint union of $k$ copies of $\mathbb R$.  Define an equivalence
relation by
\[
(x,i)\sim(x,j)
\qquad\Longleftrightarrow\qquad
x\neq 0,
\]
and let
\[
L_k
:=
\left(\bigsqcup_{i=1}^k\mathbb R_i\right)\big/\!\sim .
\]
For $x\neq0$ we denote the common equivalence class by $x$, while
\[
o_i:=[(0,i)]
\]
denotes the $i$th origin.

Equivalently,
\[
L_k
=
(\mathbb R\setminus\{0\})
\sqcup
\{o_1,\ldots,o_k\}.
\]
A neighbourhood basis at $o_i$ is given by
\[
U_i(\varepsilon)
=
(-\varepsilon,0)
\cup\{o_i\}
\cup(0,\varepsilon),
\qquad
\varepsilon>0,
\]
together with the ordinary Euclidean open sets away from the exceptional
fibre.

Set $F_k:=\{o_1,\ldots,o_k\}$; its subspace topology is discrete.

\begin{proposition}
\label{prop:Lk-basic-topology}
The space $L_k$ is a $T_1$ non-Hausdorff one-manifold.
\end{proposition}

Indeed, each origin is closed, while any two origin neighbourhoods meet in regular points arbitrarily close to $0$.

\subsection{Fundamental group of \texorpdfstring{$L_k$}{Lk}}
\label{subsec:Lk-cycles}

The non-Hausdorff exceptional fibre creates nontrivial loops even though
each individual origin chart is contractible.

For $i=1,\ldots,k$, let
\[
e_i:[0,1]\longrightarrow L_k
\]
be any path from $-1$ to $1$ whose image passes through $o_i$ and otherwise
follows the common regular line.  For example, after a reparametrization
one may take the path corresponding to the ordinary interval
$[-1,1]\subset U_i$.

The loop obtained by crossing from the negative side through $o_i$ and
returning from the positive side through $o_1$ is
\[
g_i
:=
e_i * \overline{e_1},
\qquad
i=2,\ldots,k,
\]
where $\overline{e_1}$ denotes the reverse path.

For $k=2$, the identification $\pi_1(L_2)\cong\mathbb Z$ is classical
\citep{CalcutGompfMcCarthy2012}.  The groupoid argument below gives the corresponding
$k$-origin formula.

\begin{theorem}[Fundamental group of the multiple-origin line]
\label{thm:pi1-Lk}
For $k\geq2$,
\[
\pi_1(L_k,-1)
\cong
F_{k-1},
\]
the free group on $k-1$ generators.  Under this identification the
classes
\[
[g_2],\ldots,[g_k]
\]
form a free basis.
\end{theorem}

\begin{proof}
For each $i$ set
\[
V_i
:=
(\mathbb R\setminus\{0\})\cup\{o_i\}.
\]
Then
\[
V_i\cong\mathbb R,
\]
so each $V_i$ is simply connected, and
\[
L_k=\bigcup_{i=1}^kV_i.
\]
For $i\neq j$,
\[
V_i\cap V_j
=
\mathbb R\setminus\{0\}
=
(-\infty,0)\sqcup(0,\infty).
\]

Use the basepoint set
\[
B=\{-1,1\}.
\]
Each path component of every nonempty finite intersection of the $V_i$
contains a point of $B$.  The van Kampen theorem for fundamental
groupoids therefore applies to this cover.

The groupoid $\Pi_1(V_i,B)$ has two objects, $-1$ and $1$, and, since
$V_i$ is simply connected, a unique homotopy class
\[
[e_i]:-1\longrightarrow1.
\]
By contrast,
\[
V_i\cap V_j
=
(-\infty,0)\sqcup(0,\infty)
\]
contains no path joining $-1$ to $1$.  Its fundamental groupoid on $B$
therefore contributes only the identity morphisms at the two objects.

Consequently the van Kampen pushout is the free groupoid on the $k$
parallel arrows
\[
[e_1],\ldots,[e_k]:
-1\longrightarrow1.
\]
Choosing $e_1$ as a reference arrow identifies the vertex group at
$-1$ with the free group generated by
\[
[e_i*\overline{e_1}],
\qquad
i=2,\ldots,k.
\]
Hence
\[
\pi_1(L_k,-1)\cong F_{k-1}.
\]
\end{proof}

\begin{corollary}
\label{cor:H1-Lk}
The first singular homology group is
\[
H_1(L_k;\mathbb Z)
\cong
\mathbb Z^{k-1}.
\]
\end{corollary}

\begin{proof}
This is the abelianization of
\[
\pi_1(L_k)\cong F_{k-1}.
\]
\end{proof}

\subsection{Collapse map and over-base automorphisms}
\label{subsec:Lk-collapse}

Define
\[
q_k:L_k\longrightarrow\mathbb R
\]
by
\[
q_k(x)=x
\qquad(x\neq0),
\]
and
\[
q_k(o_i)=0.
\]

\begin{proposition}
\label{prop:Lk-local-homeomorphism}
The map
\[
q_k:L_k\longrightarrow\mathbb R
\]
is a surjective local homeomorphism.  For $k>1$ it is not a covering
projection.
\end{proposition}

\begin{proof}
Away from the exceptional fibre the map is locally the identity.  At an
origin $o_i$, restriction gives
\[
q_k|_{U_i(\varepsilon)}:
U_i(\varepsilon)
\longrightarrow
(-\varepsilon,\varepsilon),
\]
which is a homeomorphism.  Thus $q_k$ is a local homeomorphism.

Suppose that some neighbourhood $U$ of $0$ were evenly covered.  Since
$q_k^{-1}(0)$ contains the $k$ distinct origins, the set
$q_k^{-1}(U)$ would have to decompose into pairwise disjoint open sheets
\[
W_1,\ldots,W_k,
\]
each mapped homeomorphically onto $U$ and with $o_i\in W_i$ after
reindexing.

For every nonzero $x\in U$, however, the fibre $q_k^{-1}(x)$ consists of
a single point.  Surjectivity of each restriction
\[
q_k|_{W_i}:W_i\to U
\]
would therefore force this same point $x$ to lie in every $W_i$,
contradicting their disjointness.  Hence no neighbourhood of $0$ is
evenly covered.
\end{proof}

Since $q_k$ is not a covering projection, it does not carry the canonical
monodromy action supplied by covering-space theory.

\begin{proposition}[Over-base automorphisms]
\label{prop:Lk-over-base-aut}
For $k\geq2$,
\[
\operatorname{Aut}_{/\mathbb R}(q_k)
:=
\{h\in\operatorname{Homeo}(L_k):q_k\circ h=q_k\}
\cong
\mathfrak S_k.
\]
\end{proposition}

\begin{proof}
Every nonzero fibre of $q_k$ is a singleton, so an over-base
homeomorphism fixes every point of
\[
\mathbb R\setminus\{0\}.
\]
It may therefore act only on
\[
F_k=\{o_1,\ldots,o_k\}.
\]
Any such action must be a permutation.

Conversely, every permutation
\[
\sigma\in\mathfrak S_k
\]
extends to a homeomorphism of $L_k$ by fixing the regular locus and
sending
\[
o_i\longmapsto o_{\sigma(i)}.
\]
The defining neighbourhood bases are permuted accordingly, and the
extension lies over $\mathbb R$.  Hence the full over-base automorphism
group is $\mathfrak S_k$.
\end{proof}

Thus $\mathfrak S_k$ is an over-base exceptional-fibre symmetry, not a covering-space monodromy group.

\section{Peripheral monodromy and component germs}
\label{sec:peripheral-monodromy}

Choose pairwise disjoint coordinate discs $U_i\ni z_i$ and a nested cofinal family
$U_i(r)$, $0<r<1$, with $U_i(r)^\times=U_i(r)\setminus\{z_i\}$.  Write
$C_i(r):=\pi_0((\pi^\circ)^{-1}(U_i(r)^\times))$.  For $0<s<r$, restriction induces
$C_i(r)\to C_i(s)$.  Connected finite covers of a
punctured disc are classified by finite-index subgroups of
$\pi_1(D^\times)\cong\mathbb Z$, so these maps are bijections
\citep{Hatcher2002,Spanier1966}.  We therefore define the \emph{peripheral component-germ set}
by
\[
C_i:=\varprojlim_{r\downarrow0} C_i(r),
\]
and identify it with $C_i(r)$ for any sufficiently small $r$.  For $c\in C_i$,
$T_{i,c}(r)$ denotes the corresponding nested tail.  Cofinality of neighbourhood bases
shows that $C_i$ and the tail filter of $c$ are independent of the chosen coordinate disc.
We also use ``peripheral end'' as shorthand for a peripheral component germ; it does
not mean a Freudenthal end.  Fix a sufficiently small reference radius $r_i^*>0$ and,
for $c\in C_i$, set $E_{i,c}:=T_{i,c}(r_i^*)$.  This is a connected representative
of the germ, containing $T_{i,c}(r)$ for all sufficiently small $r$.  All later uses of
$E_{i,c}$ refer to this fixed representative; the statements depend only on the germ.

Fix $x_0\in X^\circ$ and let $A=(\pi^\circ)^{-1}(x_0)$.  If $\gamma_i$ is a positive
peripheral loop and $\sigma_i$ its monodromy permutation, the component germs are the
cycles of $\sigma_i$.  We write
\[
\nu_i(c):=\text{the length of the cycle corresponding to }c.
\]
Equivalently, the restriction over $c$ is the connected degree-$\nu_i(c)$ cover of a
punctured disc and is locally equivalent to $w\mapsto w^{\nu_i(c)}$.  If
$H_i=\langle\gamma_i\rangle\leq G:=\pi_1(X^\circ,x_0)$, then
\[
C_i\cong H_i\backslash A,
\qquad
\nu_i(H_i\cdot a)=|H_i\cdot a|.
\]
The exceptional fibre and endpoint incidence are additional completion data beyond
ordinary peripheral monodromy.

\section{Incidence completions}
\label{sec:incidence-completions}

For each $z_i\in\Sigma$, let $C_i$ be its peripheral component-germ set and let
$F_i$ be a nonempty topological space.  The completion data specify which germs in
$C_i$ converge to which points of $F_i$.

\subsection{Closed incidence data}
\label{subsec:closed-incidence-data}

\begin{definition}[Incidence datum]
\label{def:incidence-datum}
An \emph{incidence datum} for the regular covering
\[
\pi^\circ:Y^\circ\to X^\circ
\]
consists, for every $i$ and every $c\in C_i$, of a subset
\[
A_{i,c}\subseteq F_i.
\]

The datum is called:
\begin{enumerate}[label=(\roman*)]
\item \emph{closed} if every $A_{i,c}$ is closed in $F_i$;
\item \emph{end-complete} if
\[
A_{i,c}\neq\varnothing
\qquad
\text{for every }i,c;
\]
\item \emph{endpoint-covering} if
\[
F_i=\bigcup_{c\in C_i}A_{i,c}
\qquad
\text{for every }i;
\]
\item \emph{admissible} if it is closed, end-complete, and endpoint-covering.
\end{enumerate}
\end{definition}

The corresponding incidence relation is
\[
R_i\subseteq C_i\times F_i,
\qquad
(c,a)\in R_i
\quad\Longleftrightarrow\quad
a\in A_{i,c}.
\]

For $a\in F_i$, write $I_i(a):=\{c\in C_i:a\in A_{i,c}\}$ for its
incident-end set.  If $V\subseteq F_i$ is open, set
$C_i(V):=\{c\in C_i:A_{i,c}\cap V\neq\varnothing\}$; equivalently,
$C_i(V)=\bigcup_{a\in V}I_i(a)$.

End-completeness means that every regular peripheral end has at least one
permitted exceptional endpoint.  Endpoint-covering means that every exceptional
point is approached by at least one peripheral end.

\subsection{The tail-saturated topology}
\label{subsec:tail-saturated-topology}

Set
\[
\overline Y_{\mathcal D}
:=
Y^\circ
\sqcup
\bigsqcup_{i=1}^nF_i,
\]
where
\[
\mathcal D
=
\bigl(\{F_i\},\{A_{i,c}\}\bigr)
\]
denotes the incidence datum.

For an open subset $V\subseteq F_i$ and a choice of radii
\[
\mathbf r=(r_c)_{c\in C_i(V)},
\qquad
0<r_c<1,
\]
define
\[
B_i(V,\mathbf r)
:=
V
\cup
\bigcup_{c\in C_i(V)}
T_{i,c}(r_c).
\]

\begin{definition}[Incidence completion]
\label{def:incidence-completion}
The \emph{incidence completion} associated with $\mathcal D$ is
$\overline Y_{\mathcal D}$ equipped with the topology generated by
\begin{enumerate}[label=(\roman*)]
\item all open subsets of $Y^\circ$;
\item all sets
\[
B_i(V,\mathbf r)
\]
defined above.
\end{enumerate}
The completed map
\[
\overline\pi_{\mathcal D}:
\overline Y_{\mathcal D}
\longrightarrow X
\]
is defined by
\[
\overline\pi_{\mathcal D}|_{Y^\circ}
=
\pi^\circ
\]
and
\[
\overline\pi_{\mathcal D}(F_i)=\{z_i\}.
\]
\end{definition}

\begin{proposition}[Basis theorem]
\label{prop:incidence-basis}
The family in Definition~\ref{def:incidence-completion} is a basis for a
topology on $\overline Y_{\mathcal D}$.
\end{proposition}

\begin{proof}
The family covers $Y^\circ$ by ordinary open subsets.  Every
$a\in F_i$ lies in
\[
B_i(F_i,\mathbf r)
\]
for any choice of radii.

Intersections involving an ordinary open subset of $Y^\circ$ are handled
inside the regular locus.  Incidence neighbourhoods belonging to
different marked points are disjoint after the coordinate discs
$U_i$ have been chosen pairwise disjoint.

Consider
\[
B_i(V,\mathbf r)
\qquad\text{and}\qquad
B_i(W,\mathbf s)
\]
over the same $z_i$.

If
\[
x\in
B_i(V,\mathbf r)\cap B_i(W,\mathbf s)
\cap Y^\circ,
\]
then an ordinary open neighbourhood of $x$ inside the intersection
suffices.

Suppose instead that
\[
x\in F_i.
\]
Then
\[
x\in V\cap W.
\]
Choose an open set
\[
x\in Q\subseteq V\cap W.
\]
For every
\[
c\in C_i(Q),
\]
one has
\[
c\in C_i(V)\cap C_i(W).
\]
Set
\[
t_c:=\min\{r_c,s_c\}.
\]
Then
\[
B_i(Q,\mathbf t)
\subseteq
B_i(V,\mathbf r)\cap B_i(W,\mathbf s).
\]
This verifies the basis condition.
\end{proof}

\subsection{Elementary properties}
\label{subsec:incidence-elementary}

\begin{proposition}
\label{prop:incidence-elementary}
For every incidence datum $\mathcal D$:
\begin{enumerate}[label=(\roman*)]
\item $Y^\circ$ is an open subspace of $\overline Y_{\mathcal D}$ with
      its original topology;
\item every $F_i$ is embedded with its prescribed topology;
\item
\[
\overline\pi_{\mathcal D}^{-1}(X^\circ)=Y^\circ;
\]
\item the map
\[
\overline\pi_{\mathcal D}:
\overline Y_{\mathcal D}\to X
\]
is continuous.
\end{enumerate}
If the datum is endpoint-covering, then $Y^\circ$ is dense in
$\overline Y_{\mathcal D}$.
\end{proposition}

\begin{proof}
The first assertion is built into the definition.

For the exceptional fibre,
\[
B_i(V,\mathbf r)\cap F_i=V,
\]
while neighbourhoods associated with other marked points do not meet
$F_i$.  Hence the induced subspace topology is precisely the prescribed
topology of $F_i$.

The equality
\[
\overline\pi_{\mathcal D}^{-1}(X^\circ)=Y^\circ
\]
is immediate.

To prove continuity, it is enough to work in a coordinate disc $U_i$.
For $0<r<1$,
\[
\overline\pi_{\mathcal D}^{-1}(U_i(r))
=
F_i
\cup
\bigcup_{c\in C_i}T_{i,c}(r).
\]
The tails with $A_{i,c}\neq\varnothing$ occur in the basic neighbourhood
\[
B_i(F_i,(r)_c),
\]
while any remaining tails are open subsets of the regular locus.  Hence
the inverse image is open.  Away from $\Sigma$, the map is the ordinary
covering $\pi^\circ$.

Finally suppose the datum is endpoint-covering.  If
\[
a\in F_i
\]
and
\[
B_i(V,\mathbf r)
\]
is a basic neighbourhood of $a$, endpoint-covering gives some $c\in C_i$ with
\[
a\in A_{i,c}.
\]
Since $a\in V$,
\[
A_{i,c}\cap V\neq\varnothing,
\]
so
\[
c\in C_i(V).
\]
The neighbourhood therefore contains the nonempty regular tail
$T_{i,c}(r_c)$.  Thus every neighbourhood of every exceptional point
meets $Y^\circ$.
\end{proof}

\subsection{Exact realization of endpoint sets}
\label{subsec:exact-endpoints}

\begin{definition}[Endpoint set]
\label{def:peripheral-endpoint-set}
For $c\in C_i$, define
\[
\operatorname{Ends}_{\overline Y_{\mathcal D}}(c)
\]
to be the set of all $a\in F_i$ such that the tail filter of $c$ converges
to $a$.  Explicitly,
\[
a\in\operatorname{Ends}_{\overline Y_{\mathcal D}}(c)
\]
if and only if for every neighbourhood $N$ of $a$ there exists
$r_0>0$ such that
\[
T_{i,c}(r)\subseteq N
\qquad
(0<r<r_0).
\]
\end{definition}

\begin{theorem}[Exact endpoint realization]
\label{thm:exact-endpoint-realization}
Let $\mathcal D$ be a closed incidence datum.  Then for every
$i$ and every $c\in C_i$,
\[
\operatorname{Ends}_{\overline Y_{\mathcal D}}(c)
=
A_{i,c}.
\]
\end{theorem}

\begin{proof}
Let
\[
a\in A_{i,c}.
\]
Consider a basic neighbourhood
\[
B_i(V,\mathbf r)
\]
of $a$.  Since
\[
a\in V\cap A_{i,c},
\]
one has
\[
c\in C_i(V).
\]
Therefore the neighbourhood contains
\[
T_{i,c}(r_c),
\]
and every sufficiently short $c$-tail lies in the neighbourhood.  Hence
the tail filter converges to $a$.

Conversely suppose
\[
a\notin A_{i,c}.
\]
Since $A_{i,c}$ is closed,
\[
V:=F_i\setminus A_{i,c}
\]
is an open neighbourhood of $a$.  For this $V$,
\[
A_{i,c}\cap V=\varnothing,
\]
so
\[
c\notin C_i(V).
\]
Every basic neighbourhood
\[
B_i(V,\mathbf r)
\]
therefore contains no terminal $c$-tail.  Hence the $c$-tail does not
converge to $a$.
\end{proof}

\begin{corollary}[Intrinsic incidence]
\label{cor:intrinsic-incidence}
For a closed incidence completion, the relation
\[
R_i\subseteq C_i\times F_i
\]
is recoverable from the completed topology and the fixed regular
covering:
\[
(c,a)\in R_i
\iff
a\in\operatorname{Ends}_{\overline Y_{\mathcal D}}(c).
\]
\end{corollary}

\subsection{Separation of exceptional points}
\label{subsec:incidence-separation}

\begin{proposition}[Exceptional Hausdorff separation]
\label{prop:exceptional-Hausdorff-separation}
Let $a,b\in F_i$ be distinct.  They admit disjoint neighbourhoods in
$\overline Y_{\mathcal D}$ if and only if there exist disjoint open sets
\[
a\in V\subseteq F_i,
\qquad
b\in W\subseteq F_i
\]
such that
\[
C_i(V)\cap C_i(W)=\varnothing.
\]
\end{proposition}

\begin{proof}
Suppose first that such $V$ and $W$ exist.  Choose arbitrary sufficiently
small radii on $C_i(V)$ and $C_i(W)$.  The exceptional parts $V$ and $W$
are disjoint.  Since the peripheral components are pairwise disjoint and
the two support sets of ends are disjoint, the corresponding basic
neighbourhoods
\[
B_i(V,\mathbf r)
\qquad\text{and}\qquad
B_i(W,\mathbf s)
\]
are disjoint.

Conversely suppose $a$ and $b$ have disjoint neighbourhoods.  Shrinking
them to basic neighbourhoods gives
\[
B_i(V,\mathbf r)
\cap
B_i(W,\mathbf s)
=
\varnothing.
\]
Then necessarily
\[
V\cap W=\varnothing.
\]
If some
\[
c\in C_i(V)\cap C_i(W),
\]
both basic neighbourhoods contain sufficiently short tails of the same
component $c$, and hence intersect.  Therefore
\[
C_i(V)\cap C_i(W)=\varnothing.
\]
\end{proof}

\begin{corollary}[Discrete-fibre criterion]
\label{cor:discrete-Hausdorff-separation}
Assume $F_i$ is discrete.  Distinct points $a,b\in F_i$ can be
Hausdorff-separated in the completion if and only if
\[
I_i(a)\cap I_i(b)=\varnothing.
\]
Equivalently, no peripheral end has both $a$ and $b$ as endpoints.
\end{corollary}

\section{Lifting through incidence defects}
\label{sec:incidence-lifting}

Let $f:Z\to X$ be continuous, put $D:=f^{-1}(\Sigma)$ and
$D_i:=f^{-1}(z_i)$, and fix a regular lift
$\widetilde f^\circ:Z\setminus D\to Y^\circ$ of $f|_{Z\setminus D}$.  Since
$\Sigma$ is finite and $X$ is Hausdorff, each $D_i$ is closed in $Z$ and clopen in the
subspace $D$.

\subsection{Accumulated peripheral ends}
\label{subsec:accumulated-ends}

Fix $z\in D_i$.

\begin{definition}[Accumulated-end set]
\label{def:accumulated-end-set}
The \emph{accumulated-end set} of the regular lift at $z$ is
\[
\Lambda_{\widetilde f^\circ}(z)
:=
\left\{
c\in C_i:
z\in
\overline{
(\widetilde f^\circ)^{-1}(E_{i,c})
}^{\,Z}
\right\}.
\]
\end{definition}

Thus
\[
c\in\Lambda_{\widetilde f^\circ}(z)
\]
precisely when arbitrarily close regular points of the probe are lifted
through the peripheral component $E_{i,c}$.

Because $C_i$ is finite,
\[
\Lambda_{\widetilde f^\circ}(z)
\]
is finite.

For a subset
\[
J\subseteq C_i
\]
write
\[
A_{i,J}
:=
\bigcap_{c\in J}A_{i,c},
\]
with the convention
\[
A_{i,\varnothing}:=F_i.
\]

The set
\[
A_{i,\Lambda_{\widetilde f^\circ}(z)}
\]
will be called the \emph{admissible endpoint set} at $z$.

\subsection{The incidence lift theorem}
\label{subsec:incidence-lift-theorem}

An exceptional colouring is a map
\[
\alpha:D\longrightarrow\bigsqcup_{i=1}^nF_i
\]
such that
\[
\alpha(D_i)\subseteq F_i.
\]
Since the $D_i$ are clopen in $D$, such a map is continuous if and only if
every restriction
\[
\alpha_i:=\alpha|_{D_i}:D_i\to F_i
\]
is continuous.

\begin{theorem}[Incidence lift theorem]
\label{thm:incidence-lift}
Let $\mathcal D$ be an admissible incidence datum and let
\[
\widetilde f^\circ:
Z\setminus D\to Y^\circ
\]
be a fixed regular lift of $f$.

An exceptional colouring
\[
\alpha:D\to\bigsqcup_iF_i
\]
extends $\widetilde f^\circ$ to a continuous lift
\[
\widetilde f:
Z\longrightarrow\overline Y_{\mathcal D}
\]
if and only if the following conditions hold:
\begin{enumerate}[label=(\roman*)]
\item each
\[
\alpha_i:D_i\to F_i
\]
is continuous;
\item for every $z\in D_i$,
\[
\alpha_i(z)
\in
\bigcap_{c\in
\Lambda_{\widetilde f^\circ}(z)}
A_{i,c}.
\]
\end{enumerate}
For the prescribed pair
\[
(\widetilde f^\circ,\alpha)
\]
the extension is unique.
\end{theorem}

\begin{proof}
Suppose first that
\[
\widetilde f:Z\to\overline Y_{\mathcal D}
\]
is a continuous extension.  Its restriction to $D_i$ is continuous
because $F_i$ has its prescribed subspace topology.

Fix
\[
z\in D_i
\]
and let
\[
c\in\Lambda_{\widetilde f^\circ}(z).
\]
Set
\[
a:=\widetilde f(z)\in F_i.
\]
Assume, toward a contradiction, that
\[
a\notin A_{i,c}.
\]
Since $A_{i,c}$ is closed, there exists an open neighbourhood
\[
V\subseteq F_i
\]
of $a$ such that
\[
V\cap A_{i,c}=\varnothing.
\]
Consequently
\[
c\notin C_i(V).
\]

Choose a basic neighbourhood
\[
B_i(V,\mathbf r)
\]
of $a$.  By construction it contains no terminal tail of the component
$c$.

Continuity of $\widetilde f$ at $z$ implies that some neighbourhood
$W$ of $z$ satisfies
\[
\widetilde f(W)
\subseteq
B_i(V,\mathbf r).
\]
But
\[
c\in\Lambda_{\widetilde f^\circ}(z)
\]
means that $W$ contains regular points $w$ with
\[
\widetilde f^\circ(w)\in E_{i,c}.
\]
After shrinking $W$ if necessary, continuity of $f$ forces such points
arbitrarily far down the $c$-tail.  None of these points can belong to
$B_i(V,\mathbf r)$, a contradiction.  Hence
\[
a\in A_{i,c}.
\]
Since this holds for every accumulated end $c$,
\[
a\in
\bigcap_{c\in\Lambda_{\widetilde f^\circ}(z)}
A_{i,c}.
\]

Conversely, assume (i) and (ii).  Define
\[
\widetilde f(z)
=
\begin{cases}
\widetilde f^\circ(z),
&
z\notin D,
\\[1mm]
\alpha_i(z),
&
z\in D_i.
\end{cases}
\]
The map lies over $f$ by construction.  Continuity away from $D$ is
already known.

Fix
\[
z\in D_i
\]
and let
\[
B_i(V,\mathbf r)
\]
be a basic neighbourhood of
\[
\alpha_i(z).
\]

By continuity of $\alpha_i$, there exists a neighbourhood $W_0$ of $z$
in $Z$ such that
\[
\alpha_i(W_0\cap D_i)\subseteq V.
\]
Because $\Sigma$ is finite and $f(z)=z_i$, after shrinking $W_0$ we may also
assume
\[
W_0\cap D_j=\varnothing
\qquad (j\neq i).
\]

Put
\[
C_i^-:=C_i\setminus C_i(V).
\]
If
\[
c\in C_i^-,
\]
then
\[
A_{i,c}\cap V=\varnothing.
\]
Since
\[
\alpha_i(z)\in V,
\]
condition (ii) implies
\[
c\notin\Lambda_{\widetilde f^\circ}(z).
\]
Therefore there exists a neighbourhood $W_c$ of $z$ satisfying
\[
W_c\cap
(\widetilde f^\circ)^{-1}(E_{i,c})
=
\varnothing.
\]
Because $C_i^-$ is finite, the intersection
\[
W_1
:=
W_0\cap
\bigcap_{c\in C_i^-}W_c
\]
is still a neighbourhood of $z$.

For the remaining ends $c\in C_i(V)$, when $C_i(V)\neq\varnothing$ set
\[
r_0:=\min\left\{r_i^*,\min_{c\in C_i(V)}r_c\right\}.
\]
Continuity of $f$ gives a neighbourhood $W_2$ of $z$ such that every
regular point of $W_2$ maps into
\[
U_i(r_0)^\times.
\]

Hence, for
\[
w\in
W_1\cap W_2\setminus D,
\]
the lift $\widetilde f^\circ(w)$ lies in a peripheral component belonging
to $C_i(V)$ and sufficiently far down its tail to lie in
\[
B_i(V,\mathbf r).
\]
Defect points sufficiently near $z$ are mapped by $\alpha_i$ into $V$.
Thus
\[
\widetilde f(W_1\cap W_2)
\subseteq
B_i(V,\mathbf r),
\]
proving continuity at $z$.

If
\[
C_i(V)=\varnothing,
\]
then $C_i^-=C_i$.  Intersect $W_1$ with a neighbourhood $W_2$ on which
\[
f(W_2)\subseteq U_i(r_i^*).
\]
Every regular point of $W_2$ would then lift into one of the representatives
$E_{i,c}$, while $W_1$ excludes all such representatives.  Hence
$W_1\cap W_2$ contains no regular points, and continuity reduces to that of
$\alpha_i$.

The extension is unique once its values on $Z\setminus D$ and on $D$ are
prescribed.
\end{proof}

\subsection{Higher incidence}
\label{subsec:higher-incidence-lifting}

Suppose that a defect point is approached through several peripheral
components:
\[
\Lambda_{\widetilde f^\circ}(z)=J
\subseteq C_i.
\]
Then the admissible exceptional values form $A_{i,J}=\bigcap_{c\in J}A_{i,c}$.

\begin{definition}[Incidence nerve]
\label{def:incidence-nerve}
The \emph{incidence nerve} at $z_i$ is
\[
N_i(\mathcal D)
:=
\left\{
J\subseteq C_i:
\bigcap_{c\in J}A_{i,c}\neq\varnothing
\right\}.
\]
\end{definition}

\begin{proposition}
\label{prop:incidence-nerve-complex}
The collection
\[
N_i(\mathcal D)
\]
is an abstract simplicial complex on vertex set $C_i$.
\end{proposition}

\begin{proof}
If
\[
J\in N_i(\mathcal D)
\]
and
\[
K\subseteq J,
\]
then
\[
\bigcap_{c\in J}A_{i,c}
\subseteq
\bigcap_{c\in K}A_{i,c}.
\]
Hence the latter intersection is nonempty.
\end{proof}

The corresponding intersection diagram is
\[
\mathcal A_i:
\mathcal P(C_i)^{\mathrm{op}}
\longrightarrow
\mathbf{Top},
\qquad
J\longmapsto A_{i,J}.
\]

For
\[
J\subseteq K,
\]
there is a canonical inclusion
\[
A_{i,K}\hookrightarrow A_{i,J}.
\]

\subsection{Finite discrete lift counts}
\label{subsec:discrete-lift-counts}

Assume now that $F_i$ is finite and discrete.

Define
\[
\lambda_i(J)
:=
\left|
A_{i,J}
\right|
=
\left|
\bigcap_{c\in J}A_{i,c}
\right|.
\]

\begin{corollary}[Local lift count]
\label{cor:local-lift-count}
Let
\[
z\in D_i
\]
be an isolated defect point of the probe.  For a fixed regular lift
$\widetilde f^\circ$, the number of possible completed values at $z$ is
\[
\lambda_i
\bigl(
\Lambda_{\widetilde f^\circ}(z)
\bigr).
\]
\end{corollary}

If the defect locus consists of finitely many isolated points
\[
D=\{z_1,\ldots,z_m\},
\]
then the choices are independent.

\begin{corollary}[Finite defect counting]
\label{cor:finite-defect-counting}
For a fixed regular lift and finite discrete exceptional fibres,
\[
\#\operatorname{Ext}(\widetilde f^\circ)
=
\prod_{j=1}^m
\left|
\bigcap_{c\in
\Lambda_{\widetilde f^\circ}(z_j)}
A_{i(j),c}
\right|,
\]
where
\[
f(z_j)=z_{i(j)}.
\]
\end{corollary}

\section{Discrete incidence graphs and local topology}
\label{sec:discrete-incidence-topology}

Assume in this section that the exceptional fibre is finite and discrete.  Then
the endpoint relation is a finite bipartite graph, and the topology of a
completed neighbourhood can be computed directly from that graph.

Fix one marked point $z_i$ and suppress the index $i$.  Let
$C=\{c_1,\ldots,c_s\}$ be the peripheral-end set, let
$\nu:C\to\mathbb N_{>0}$ be the degree function, and let
$F=\{a_1,\ldots,a_k\}$ be a finite discrete exceptional fibre with admissible endpoint
sets $A_c\subseteq F$.  Thus each $A_c$ is nonempty and
$F=\bigcup_{c\in C}A_c$.

Choose a sufficiently small punctured coordinate disc and a common
radius $0<r<1$.  Write
\[
E_c:=T_c(r)
\]
for the corresponding peripheral tails and define the local completed
neighbourhood
\[
\mathcal N
:=
F\cup\bigcup_{c\in C}E_c.
\]
Shrinking $r$ changes neither the incidence relation nor the homotopy
calculations below.

\subsection{The incidence graph}
\label{subsec:incidence-graph}

\begin{definition}[Incidence graph]
\label{def:incidence-graph}
The \emph{incidence graph} is the finite bipartite graph
\[
G_{\mathrm{inc}}
=
(C\sqcup F,E),
\]
where
\[
(c,a)\in E
\quad\Longleftrightarrow\quad
a\in A_c.
\]
The vertex $c\in C$ carries the additional label
\[
\nu(c).
\]
The resulting pair
\[
(G_{\mathrm{inc}},\nu)
\]
will be called the \emph{monodromy-labelled incidence graph}.
\end{definition}

Admissibility says exactly that $G_{\mathrm{inc}}$ has no isolated
vertices.

The bipartition records different types of data: $c\in C$ denotes a
peripheral component germ and $a\in F$ an exceptional point, with
$c---a$ exactly when the $c$-tail converges to $a$.  Set
\[
d_c:=|A_c|=\deg_G(c),
\qquad
\deg_G(a)=|\{c\in C:a\in A_c\}|.
\]
Thus $d_c>1$ records endpoint multiplicity for one germ, while
$\deg_G(a)>1$ records incidence of several germs at one exceptional point.

\subsection{A canonical open cover}
\label{subsec:incidence-open-cover}

For every $a\in F$, define
\[
U_a
:=
\{a\}
\cup
\bigcup_{\substack{c\in C\\a\in A_c}}E_c.
\]

\begin{lemma}
\label{lem:Ua-open-contractible}
Each $U_a$ is an open contractible subspace of $\mathcal N$.
\end{lemma}

\begin{proof}
Since $F$ is discrete, $\{a\}$ is open in the exceptional fibre.
A basic incidence neighbourhood of $a$ contains sufficiently short
tails of precisely the ends incident with $a$, and hence is contained in
$U_a$.  Every regular point of $U_a$ has an ordinary open neighbourhood
inside its peripheral component.  Thus $U_a$ is open.

For each incident end $c$, the subspace
\[
E_c\cup\{a\}
\]
is a disc: in the cyclic local coordinate it is obtained from a punctured
disc by adjoining one apex.  The different such discs meet only in the
common point $a$.

Choose radial coordinates independently on the incident peripheral
components.  Simultaneously contracting each radial coordinate to the
common apex gives a contraction of $U_a$ to $a$.
\end{proof}

The sets $U_a$ cover $\mathcal N$.  Indeed every exceptional point lies
in its own $U_a$, and every peripheral end $c$ has at least one endpoint
by end-completeness.

\begin{lemma}
\label{lem:Ua-intersections}
Let
\[
S=\{a_0,\ldots,a_p\}\subseteq F
\]
contain distinct exceptional points and suppose $p\geq1$.  Then
\[
\bigcap_{a\in S}U_a
=
\bigsqcup_{\substack{c\in C\\S\subseteq A_c}}E_c.
\]
Consequently every connected component of a nonempty intersection of at
least two members of the cover is homotopy equivalent to $S^1$.
\end{lemma}

\begin{proof}
No exceptional point belongs to two distinct sets $U_a$ and $U_b$.
A regular point belongs to $U_a$ precisely when it lies in a peripheral
component $E_c$ satisfying
\[
a\in A_c.
\]
Thus it belongs to all $U_a$ with $a\in S$ precisely when
\[
S\subseteq A_c.
\]
Distinct peripheral components are disjoint.  Finally,
\[
E_c\cong D^\times\simeq S^1.
\]
\end{proof}

\subsection{The Mayer--Vietoris spectral sequence}
\label{subsec:incidence-MVSS}

We use the homological Mayer--Vietoris spectral sequence associated with
a finite open cover \citep{McCleary2001}.  For
\[
\mathcal U=\{U_a\}_{a\in F},
\]
its first page is
\[
E^1_{p,q}
=
\bigoplus_{\substack{S\subseteq F\\|S|=p+1}}
H_q(U_S;\mathbb Z),
\qquad
U_S:=\bigcap_{a\in S}U_a,
\]
with the usual alternating Čech differential in the $p$-direction, and
it converges to
\[
H_{p+q}(\mathcal N;\mathbb Z).
\]

By Lemmas~\ref{lem:Ua-open-contractible}
and~\ref{lem:Ua-intersections},
only the rows
\[
q=0
\qquad\text{and}\qquad
q=1
\]
can be nonzero.

\subsection{The zeroth row and the incidence graph}
\label{subsec:MV-zero-row}

For $p=0$, $E^1_{0,0}\cong\mathbb Z[F]$.  For $p\geq1$, a basis element of
$E^1_{p,0}$ is a pair $(c,S)$ with $S\subseteq A_c$ and $|S|=p+1$.
Define a finite $\Delta$-complex $K_{\mathrm{inc}}$ with $0$-skeleton $F$ by
attaching, for each $c\in C$, one simplex $\Delta_c$ with vertex set $A_c$.
Positive-dimensional faces belonging to different $c$'s are retained
separately, while vertices carrying the same label in $F$ are identified.
Its cellular chain complex is the $q=0$ row of the Mayer--Vietoris spectral
sequence.

Let $G_{\mathrm{inc}}^{\mathrm{red}}$ be obtained from $G_{\mathrm{inc}}$ by
contracting the unique edge incident with every vertex $c$ satisfying
$|A_c|=1$.  For $|A_c|\geq2$, the barycentric subdivision of $\Delta_c$
deformation retracts relative to its vertices onto the star joining its
barycentre to the vertices of $A_c$.  The retractions patch because distinct
$\Delta_c$ meet only in the common $0$-skeleton.  The simplices with
$|A_c|=1$ already consist of their sole exceptional vertex, so the resulting
graph is $G_{\mathrm{inc}}^{\mathrm{red}}$.  Hence
\[
K_{\mathrm{inc}}\simeq G_{\mathrm{inc}}^{\mathrm{red}}.
\]
Each contraction defining $G_{\mathrm{inc}}^{\mathrm{red}}$ removes a leaf
edge and is a homotopy equivalence, so
\[
G_{\mathrm{inc}}^{\mathrm{red}}\simeq G_{\mathrm{inc}}.
\]
Therefore
\[
E^2_{p,0}\cong H_p(G_{\mathrm{inc}};\mathbb Z).
\]
In particular $E^2_{p,0}=0$ for $p\geq2$.  If
$r=b_0(G_{\mathrm{inc}})$, then
\[
H_0(G_{\mathrm{inc}};\mathbb Z)\cong\mathbb Z^r,
\qquad
H_1(G_{\mathrm{inc}};\mathbb Z)\cong
\mathbb Z^{|E|-|C|-|F|+r}.
\]

\subsection{The first row and endpoint multiplicity}
\label{subsec:MV-first-row}

Since every $U_a$ is contractible,
\[
E^1_{0,1}=0.
\]

For $p\geq1$, each component $E_c$ of an intersection contributes
\[
H_1(E_c;\mathbb Z)\cong\mathbb Z.
\]
The $q=1$ row therefore decomposes as a direct sum over peripheral ends:
\[
E^1_{\bullet,1}
=
\bigoplus_{c\in C}K_\bullet(c),
\]
where
\[
K_p(c)
=
\mathbb Z
\left\langle
S\subseteq A_c:\ |S|=p+1
\right\rangle
\qquad(p\geq1),
\]
and
\[
K_0(c)=0.
\]

For fixed $c$, this is the simplicial chain complex of the simplex on
$A_c$, truncated by deleting degree zero.  Since a simplex is
contractible, its only homology after this truncation occurs in degree
one and is naturally
\[
\widetilde H_0(A_c;\mathbb Z).
\]

Consequently
\[
E^2_{1,1}
\cong
\bigoplus_{c\in C}
\widetilde H_0(A_c;\mathbb Z)
\]
and
\[
E^2_{p,1}=0
\qquad
(p\neq1).
\]

For a discrete set of cardinality $d_c$,
\[
\widetilde H_0(A_c;\mathbb Z)
\cong
\mathbb Z^{d_c-1}.
\]

Thus
\[
E^2_{1,1}
\cong
\mathbb Z^{
\sum_{c\in C}(d_c-1)
}.
\]

Since
\[
\sum_{c\in C}d_c=|E|,
\]
the rank is
\[
|E|-|C|.
\]

\subsection{Collapse of the spectral sequence}
\label{subsec:MV-collapse}

The only nonzero terms on the second page occur at
\[
(0,0),\qquad
(1,0),\qquad
(1,1).
\]
There are therefore no possible nonzero higher differentials.

Hence
\[
E^2=E^\infty.
\]

Because all groups involved are free abelian, no extension problem
remains.

\begin{theorem}[Discrete incidence homology]
\label{thm:discrete-incidence-homology}
Let $\mathcal N$ be a local incidence completion with finite discrete
exceptional fibre and incidence graph
\[
G_{\mathrm{inc}}=(C\sqcup F,E).
\]
Then
\[
H_0(\mathcal N;\mathbb Z)
\cong
H_0(G_{\mathrm{inc}};\mathbb Z),
\]
\[
H_1(\mathcal N;\mathbb Z)
\cong
H_1(G_{\mathrm{inc}};\mathbb Z),
\]
and
\[
H_2(\mathcal N;\mathbb Z)
\cong
\bigoplus_{c\in C}
\widetilde H_0(A_c;\mathbb Z).
\]
Moreover,
\[
H_q(\mathcal N;\mathbb Z)=0
\qquad(q\geq3).
\]

If $r=b_0(G_{\mathrm{inc}})$, then
\[
H_1(\mathcal N;\mathbb Z)
\cong
\mathbb Z^{\,|E|-|C|-|F|+r},
\]
and
\[
H_2(\mathcal N;\mathbb Z)
\cong
\mathbb Z^{\,|E|-|C|}.
\]
\end{theorem}

\subsection{Fundamental group}
\label{subsec:incidence-pi1}

\begin{theorem}[Fundamental group of a discrete incidence completion]
\label{thm:incidence-fundamental-group}
Let $G_\lambda$ be a connected component of $G_{\mathrm{inc}}$, and let
$\mathcal N_\lambda$ be the corresponding connected component of the local completion.  Then
\[
\pi_1(\mathcal N_\lambda)\cong F_{\beta_1(G_\lambda)}.
\]
In particular, if $G_{\mathrm{inc}}$ is connected,
\[
\pi_1(\mathcal N)\cong F_{|E|-|C|-|F|+1}.
\]
\end{theorem}

\begin{proof}
Use the cover $\{U_a\}_{a\in F}$ and the groupoid form of van Kampen
\citep{Brown2006}.  For each $c\in C$, choose a basepoint $x_c\in E_c$ and set
\[
B:=\{x_c:c\in C\}.
\]
Then $B\cap U_a=B_a:=\{x_c:a\in A_c\}$, and every path component of every
nonempty finite intersection of members of the cover contains the corresponding
$x_c$.  The groupoid van Kampen theorem therefore applies with the single global
object set $B$.  Since $U_a$ is contractible, $\Pi_1(U_a,B_a)$ is the pair groupoid
on $B_a$; equivalently it is the fundamental groupoid of a star with centre labelled
$a$ and leaves labelled by the incident germs $c$.

For $a\neq b$, each component of $U_a\cap U_b$ is some $E_c$ with
$a,b\in A_c$.  With the single object $x_c$ on that component its fundamental group is
$\pi_1(E_c,x_c)\cong\mathbb Z$.  The generator maps trivially into both
$\Pi_1(U_a,B_a)$ and $\Pi_1(U_b,B_b)$ because $E_c\cup\{a\}$ and
$E_c\cup\{b\}$ cap the peripheral loop.  Higher intersections introduce the same
objects $x_c$ and no additional arrows after these peripheral loops are killed.

The van Kampen colimit is therefore the free groupoid obtained by gluing, for each
$a\in F$, the star on the incident germs.  This is the fundamental groupoid of the
bipartite incidence graph $G_{\mathrm{inc}}$.  The vertex group of a connected graph is
free of rank its cycle number $\beta_1$, giving the result.
\end{proof}

\begin{corollary}
\label{cor:incidence-simply-connected}
A connected discrete incidence completion is simply connected if and only if
$G_{\mathrm{inc}}$ is a tree.  More generally, $H_1(\mathcal N)=0$ if and only if
$G_{\mathrm{inc}}$ is a forest.
\end{corollary}

\subsection{Endpoint uniqueness and second homology}
\label{subsec:H2-unique-endpoints}

\begin{corollary}
\label{cor:H2-unique-completion}
For a finite discrete exceptional fibre,
\[
H_2(\mathcal N;\mathbb Z)=0
\quad\Longleftrightarrow\quad
|A_c|=1
\text{ for every }c\in C.
\]
\end{corollary}

\begin{proof}
By Theorem~\ref{thm:discrete-incidence-homology},
\[
H_2(\mathcal N)
\cong
\bigoplus_{c\in C}
\mathbb Z^{|A_c|-1}.
\]
The group vanishes exactly when every summand vanishes.
\end{proof}

If every $c$ has a unique endpoint, every connected component of $G_{\mathrm{inc}}$ is a star; hence both $H_1$ and $H_2$ vanish.  Endpoint multiplicity is detected by $H_2$, while graph cycles are detected by $H_1$.

\subsection{Euler characteristic}
\label{subsec:incidence-euler}

\begin{corollary}
\label{cor:incidence-euler}
For every finite discrete local incidence completion,
\[
\chi(\mathcal N)=|F|.
\]
\end{corollary}

\begin{proof}
Let
\[
r=b_0(G_{\mathrm{inc}}).
\]
Then
\[
b_0=r,
\]
\[
b_1=|E|-|C|-|F|+r,
\]
and
\[
b_2=|E|-|C|.
\]
Hence
\[
\chi(\mathcal N)
=
r-
\bigl(|E|-|C|-|F|+r\bigr)
+
\bigl(|E|-|C|\bigr)
=
|F|.
\]
\end{proof}

\section{Reconstruction from lift probes}
\label{sec:lift-reconstruction}

Work locally over one marked point and suppress its index.  Let
$\pi^\circ:Y^\circ\to D^\times$ be the regular finite covering, let $C$ be its
peripheral-end set, and let $F$ be a finite discrete exceptional fibre with endpoint
sets $A_c\subseteq F$.  Write $\overline\pi:\overline Y\to D$ for the completed local
map.

\subsection{Power probes}
\label{subsec:power-probes}

For every integer $n\geq1$, define
\[
f_n:D\longrightarrow D,
\qquad
f_n(w)=w^n.
\]
Its restriction
\[
f_n^\circ:D^\times\longrightarrow D^\times
\]
induces
\[
(f_n^\circ)_*:
\pi_1(D^\times)
\longrightarrow
\pi_1(D^\times),
\qquad
1\longmapsto n.
\]

Fix a peripheral component
\[
c\in C
\]
and put
\[
m:=\nu(c).
\]
The connected covering
\[
E_c\longrightarrow D^\times
\]
corresponds to the subgroup
\[
m\mathbb Z
\leq
\mathbb Z.
\]

\begin{proposition}[Power-lift divisibility]
\label{prop:power-lift-divisibility}
The punctured probe
\[
f_n^\circ:D^\times\to D^\times
\]
admits a lift
\[
\widetilde f_n^\circ:
D^\times\longrightarrow E_c
\]
if and only if
\[
m\mid n.
\]
If a point of the fibre over a chosen basepoint is prescribed, the lift,
when it exists, is unique.
\end{proposition}

\begin{proof}
The usual covering-space lifting criterion gives a lift precisely when
\[
(f_n^\circ)_*
\pi_1(D^\times)
\subseteq
(\pi^\circ|_{E_c})_*
\pi_1(E_c).
\]
Under the standard identifications with $\mathbb Z$, this condition is
\[
n\mathbb Z
\subseteq
m\mathbb Z,
\]
which is equivalent to
\[
m\mid n.
\]
Uniqueness after fixing one value follows from ordinary unique path
lifting on the regular covering.
\end{proof}

\[
\nu(c)
=
\min
\left\{
n\geq1:
f_n^\circ
\text{ lifts to }E_c
\right\}.
\]

\subsection{Completion of a based power lift}
\label{subsec:based-power-lifts}

Suppose a basepoint
\[
w_0\in D^\times
\]
and a point
\[
y_0\in E_c
\cap
(\pi^\circ)^{-1}(f_n(w_0))
\]
have been fixed.

When
\[
\nu(c)\mid n,
\]
Proposition~\ref{prop:power-lift-divisibility} produces a unique regular
lift
\[
\widetilde f_{n,c,y_0}^\circ:D^\times\to E_c.
\]
As
\[
w\to0,
\]
this lift remains entirely in the peripheral component $c$.  Hence its
accumulated-end set at the centre is
\[
\Lambda(0)=\{c\}.
\]

\begin{theorem}[Based power-probe count]
\label{thm:based-power-probe-count}
For a fixed peripheral component $c$ and fixed starting sheet,
the number of completed lifts of
\[
f_n:D\to D
\]
through $c$ is
\[
L_c^{\mathrm{based}}(n)
=
\begin{cases}
|A_c|,
&
\nu(c)\mid n,
\\[1mm]
0,
&
\nu(c)\nmid n.
\end{cases}
\]
\end{theorem}

\begin{proof}
If
\[
\nu(c)\nmid n,
\]
there is no regular lift into $E_c$.

If
\[
\nu(c)\mid n,
\]
the chosen starting sheet determines a unique regular lift.  Its
accumulated-end set at $0$ is $\{c\}$.  Since $F$ is discrete, every
choice
\[
a\in A_c
\]
gives a continuous exceptional colouring of the singleton defect locus,
and the incidence lift theorem shows that these are precisely the
possible completions.
\end{proof}

\begin{corollary}[Recovery of one labelled end]
\label{cor:recover-labelled-end}
For a known peripheral component $c$,
\[
\nu(c)
=
\min
\left\{
n\geq1:
L_c^{\mathrm{based}}(n)>0
\right\},
\]
and
\[
|A_c|
=
L_c^{\mathrm{based}}\bigl(\nu(c)\bigr).
\]
\end{corollary}

\subsection{Multi-arm probes}
\label{subsec:multi-arm-probes}

Power probes approach the defect through one peripheral component.
To recover the intersections
\[
A_{c_1}\cap\cdots\cap A_{c_r},
\]
we use probes with several regular arms meeting at one defect point.

Let
\[
J=\{c_1,\ldots,c_q\}
\subseteq C
\]
be nonempty.

Define the $q$-armed star
\[
S_J
:=
\left(
\bigsqcup_{j=1}^q[0,1]_j
\right)\big/\!\sim,
\]
where all points $0_j$ are identified to a single central vertex
\[
v.
\]
Write
\[
(0,1]_j
\]
for the $j$th open arm.

Choose a radial path
\[
\gamma_j:[0,1]\to D
\]
on each arm with
\[
\gamma_j(0)=0
\]
and
\[
\gamma_j(t)\neq0
\qquad(t>0).
\]
Define
\[
f_J:S_J\longrightarrow D
\]
by
\[
f_J|_{[0,1]_j}=\gamma_j.
\]

For each $c_j\in J$, choose a regular lift of the $j$th punctured arm
through $E_{c_j}$.  Since the punctured arms are disjoint in the domain,
these choices define a regular lift
\[
\widetilde f_J^\circ:
S_J\setminus\{v\}
\longrightarrow
Y^\circ.
\]

By construction,
\[
\Lambda_{\widetilde f_J^\circ}(v)=J.
\]

\begin{theorem}[Multi-arm lift count]
\label{thm:multi-arm-lift-count}
For finite discrete $F$, the number of completions of the fixed regular
multi-arm lift at the central defect point is
\[
\lambda(J)
:=
\left|
\bigcap_{c\in J}A_c
\right|.
\]
\end{theorem}

\begin{proof}
The defect locus of $f_J$ is the singleton $\{v\}$ and its accumulated-end
set is exactly $J$.  The possible values at $v$ are
\[
A_J
=
\bigcap_{c\in J}A_c.
\]
Since $F$ is discrete, every choice of one point of this intersection is
a continuous colouring of the singleton defect locus.
\end{proof}

The incidence nerve introduced in
Section~\ref{sec:incidence-lifting} is the support of this counting function:
\[
J\in N(\mathcal D)
\quad\Longleftrightarrow\quad
\lambda(J)>0.
\]

\subsection{Reconstruction of the incidence graph}
\label{subsec:graph-reconstruction}

For an exceptional point
\[
a\in F,
\]
define its end-neighbourhood
\[
N(a)
:=
\{c\in C:a\in A_c\}.
\]
Endpoint-covering implies
\[
N(a)\neq\varnothing.
\]

For every nonempty subset
\[
T\subseteq C,
\]
let
\[
n_T
:=
\#\{a\in F:N(a)=T\}.
\]

The numbers
\[
n_T
\]
determine the incidence graph up to permutation of vertices in $F$ with
the same neighbourhood.

For every nonempty
\[
J\subseteq C,
\]
one has
\[
a\in
\bigcap_{c\in J}A_c
\]
if and only if
\[
J\subseteq N(a).
\]
Therefore
\[
\lambda(J)
=
\sum_{T\supseteq J}n_T.
\]

This is the zeta transform on the Boolean lattice of subsets of $C$.

\begin{theorem}[Incidence reconstruction]
\label{thm:incidence-reconstruction}
For every nonempty subset $T\subseteq C$,
\[
n_T
=
\sum_{J\supseteq T}
(-1)^{|J|-|T|}
\lambda(J).
\]
Consequently the complete collection of multi-arm lift counts
\[
\{\lambda(J):\varnothing\neq J\subseteq C\}
\]
determines the discrete incidence graph
\[
G_{\mathrm{inc}}
\]
up to relabelling of the exceptional vertices.
\end{theorem}

\begin{proof}
The relation
\[
\lambda(J)
=
\sum_{T\supseteq J}n_T
\]
is the zeta transform on the finite Boolean lattice
\[
\mathcal P(C).
\]
Its M\"obius function is
\[
\mu(J,T)
=
(-1)^{|T|-|J|}
\qquad(J\subseteq T).
\]
Boolean M\"obius inversion therefore gives
\[
n_T
=
\sum_{J\supseteq T}
(-1)^{|J|-|T|}
\lambda(J).
\]

Knowing every $n_T$ determines how many exceptional vertices have each
possible neighbourhood $T\subseteq C$.  This determines the bipartite
incidence graph up to relabelling of vertices having identical
neighbourhoods.
\end{proof}

\subsection{Combined reconstruction of the labelled local model}
\label{subsec:combined-reconstruction}

Based power probes determine
\[
\nu(c)
\qquad\text{and}\qquad
|A_c|
\]
for each labelled peripheral component.

Multi-arm probes determine
\[
\lambda(J)
=
|A_J|
\]
for every nonempty
\[
J\subseteq C,
\]
and therefore reconstruct
\[
G_{\mathrm{inc}}.
\]

\begin{theorem}[Lift reconstruction of a discrete local completion]
\label{thm:full-lift-reconstruction}
Let
\[
\overline\pi:\overline Y\to D
\]
be a finite discrete incidence completion of a finite covering of
$D^\times$.  Assume the peripheral components are labelled.

Then the collection of:
\begin{enumerate}[label=(\roman*)]
\item based power-probe lift counts
\[
L_c^{\mathrm{based}}(n),
\qquad
c\in C,\ n\geq1,
\]
and
\item multi-arm completion counts
\[
\lambda(J),
\qquad
\varnothing\neq J\subseteq C,
\]
determines the monodromy-labelled incidence graph
\[
(G_{\mathrm{inc}},\nu)
\]
up to relabelling of exceptional vertices with identical incidence
neighbourhoods.
\end{enumerate}
\end{theorem}

\section{Finite exceptional fibres}
\label{sec:finite-fibres}

Let $F$ be a finite exceptional fibre.  We use the standard equivalence between
finite topologies and specialization preorders, together with the standard description
of the Kolmogorov quotient \citep{Alexandroff1937,Stong1966,McCord1966,Barmak2011}.
Write $x\preceq y$ when $x\in\overline{\{y\}}$, let $\rho:F\to P:=F/{\sim}$ be
the Kolmogorov quotient, and put $\mu(p):=|\rho^{-1}(p)|$.  Then $P$ is a finite
poset, and every closed subset of $F$ is the inverse image of a lower set in $P$.

For an admissible endpoint family $\{A_c\}_{c\in C}$ define
\[
I(p):=\{c\in C:\rho^{-1}(p)\subseteq A_c\}.
\]
Because the $A_c$ are closed, $I$ is antitone:
\[
p\leq q\Longrightarrow I(p)\supseteq I(q).
\]
Endpoint-covering gives $I(p)\neq\varnothing$ for every $p$, and end-completeness gives
$\bigcup_{p\in P}I(p)=C$.

\begin{definition}[Finite incidence triple]
A finite incidence triple is a triple $(P,\mu,I)$ consisting of a finite poset
$P$, a weight $\mu:P\to\mathbb N_{>0}$, and an antitone map
\[
I:P^{\mathrm{op}}\longrightarrow\mathcal P(C)\setminus\{\varnothing\}
\]
with $\bigcup_p I(p)=C$.
\end{definition}

\begin{theorem}[Finite incidence classification]
\label{thm:finite-incidence-classification}
For fixed labelled peripheral-end set $C$, admissible finite incidence data
$(F,\{A_c\}_{c\in C})$ are classified up to isomorphism by finite incidence
triples $(P,\mu,I)$ up to weight- and incidence-preserving poset isomorphism.
\end{theorem}

\begin{proof}
The forward construction is the one above.  Conversely, replace each $p\in P$
by a set $F_p$ of $\mu(p)$ indistinguishable points and define the
specialization preorder by $x\preceq y$ exactly when their classes satisfy
$p_x\leq p_y$.  Set
\[
A_c:=\bigsqcup_{\{p\in P:\,c\in I(p)\}}F_p.
\]
Antitonicity makes each $A_c$ a lower set, hence closed; the two incidence
conditions give end-completeness and endpoint-covering.  The two constructions are
inverse up to the stated isomorphisms.
\end{proof}

For $J\subseteq C$ put $P_J=\{p\in P:J\subseteq I(p)\}$.  Then
\[
\bigcap_{c\in J}A_c=\bigsqcup_{p\in P_J}F_p,
\qquad
\left|\bigcap_{c\in J}A_c\right|
=\sum_{J\subseteq I(p)}\mu(p).
\]
The triple $(P,\mu,I)$ records separation order, indistinguishability multiplicity, and peripheral incidence.  The discrete incidence graph is the special case in which $P$ is an antichain and $\mu\equiv1$.

\section{Structured monodromy classification}
\label{sec:structured-monodromy}

Let $G:=\pi_1(X^\circ,x_0)$ and let $A:=(\pi^\circ)^{-1}(x_0)$ be the
regular fibre.  Path lifting gives the usual left monodromy action
$\rho:G\to\Sym(A)$.

Up to the usual choice of basepoint and labelling, the regular covering is encoded by the finite $G$-set $A$.  The peripheral subgroup $H_i$ depends on a choice of path from $x_0$ to a small loop about $z_i$ and is therefore defined up to conjugacy in $G$.  If $H_i'=gH_ig^{-1}$, the map
\[
H_i\!\cdot a
\longmapsto
H_i'\!\cdot(g\cdot a)
\]
is a well-defined bijection $H_i\backslash A\to H_i'\backslash A$ preserving orbit cardinalities and transporting the endpoint-incidence relation.  The resulting structured category is therefore independent, up to equivalence, of the chosen peripheral representatives.

\subsection{Peripheral topological incidence}
\label{subsec:peripheral-topological-incidence}

For each marked point $z_i$, let
\[
F_i
\]
be the exceptional topological fibre and let
\[
R_i\subseteq C_i\times F_i
\]
be the incidence relation
\[
(c,a)\in R_i
\quad\Longleftrightarrow\quad
a\in A_{i,c}.
\]

Using
\[
C_i\cong H_i\backslash A,
\]
we may regard $R_i$ instead as a relation
\[
R_i
\subseteq
(H_i\backslash A)\times F_i.
\]

Thus the local structured datum is $(H_i\backslash A,\nu_i;F_i,\tau_i;R_i)$, where the first two terms are determined by ordinary monodromy and the latter terms record the exceptional topology and endpoint incidence.

\begin{definition}[Structured monodromy datum]
\label{def:structured-monodromy-datum}
A \emph{structured monodromy datum} over $(X,\Sigma,x_0)$ consists of:
\begin{enumerate}[label=(\roman*)]
\item a finite $G$-set $A$;
\item for every $z_i\in\Sigma$, a nonempty topological space $F_i$;
\item for every $i$, an admissible closed incidence relation
\[
R_i\subseteq(H_i\backslash A)\times F_i,
\]
meaning that the associated endpoint sets
\[
A_{i,c}
:=
\{a\in F_i:(c,a)\in R_i\}
\]
are nonempty and closed and satisfy
\[
F_i=\bigcup_{c\in H_i\backslash A}A_{i,c}.
\]
\end{enumerate}
\end{definition}

For finite exceptional fibres, one may replace $(F_i,\tau_i,R_i)$ by the
equivalent incidence triple
\[
(P_i,\mu_i,\mathcal I_i),
\]
where
\[
\mathcal I_i:
P_i^{\mathrm{op}}
\longrightarrow
\mathcal P(H_i\backslash A)\setminus\{\varnothing\}.
\]

\subsection{Morphisms of structured data}
\label{subsec:structured-morphisms}

Let
\[
\mathfrak A
=
\left(
A;\{F_i,R_i\}_{i=1}^n
\right)
\]
and
\[
\mathfrak B
=
\left(
B;\{F_i',R_i'\}_{i=1}^n
\right)
\]
be structured monodromy data for the same base and peripheral
subgroups.

A $G$-equivariant map
\[
u:A\longrightarrow B
\]
induces, for every $i$, a map on peripheral orbit sets
\[
\overline u_i:
H_i\backslash A
\longrightarrow
H_i\backslash B,
\qquad
H_i\cdot a
\longmapsto
H_i\cdot u(a).
\]

\begin{definition}[Structured morphism]
\label{def:structured-morphism}
A \emph{structured morphism}
\[
(u,\phi_i):
\mathfrak A\longrightarrow\mathfrak B
\]
consists of:
\begin{enumerate}[label=(\roman*)]
\item a $G$-equivariant map
\[
u:A\to B;
\]
\item continuous maps
\[
\phi_i:F_i\to F_i';
\]
\item the incidence condition
\[
\phi_i(A_{i,c})
\subseteq
A_{i,\overline u_i(c)}'
\]
for every $i$ and every
\[
c\in H_i\backslash A.
\]
\end{enumerate}
\end{definition}

Equivalently,
\[
(c,a)\in R_i
\quad\Longrightarrow\quad
(\overline u_i(c),\phi_i(a))\in R_i'.
\]

\begin{proposition}
\label{prop:structured-category}
Structured monodromy data with the morphisms of
Definition~\ref{def:structured-morphism} form a category.
\end{proposition}

\begin{proof}
The identity $G$-map and identity fibre maps preserve incidence.

Suppose
\[
(u,\phi_i):
\mathfrak A\to\mathfrak B
\]
and
\[
(v,\psi_i):
\mathfrak B\to\mathfrak C
\]
are structured morphisms.  Then
\[
v\circ u
\]
is $G$-equivariant and
\[
\overline{(v\circ u)}_i
=
\overline v_i\circ\overline u_i.
\]
Moreover,
\[
\phi_i(A_{i,c})
\subseteq
A_{i,\overline u_i(c)}'
\]
and
\[
\psi_i
\left(
A_{i,\overline u_i(c)}'
\right)
\subseteq
A_{i,\overline v_i(\overline u_i(c))}''.
\]
Hence
\[
(\psi_i\circ\phi_i)(A_{i,c})
\subseteq
A_{i,\overline{(v\circ u)}_i(c)}''.
\]
Thus the composite is again structured.
\end{proof}

Denote this category by
\[
\mathbf{StrMon}(X,\Sigma,x_0).
\]

\subsection{The geometric category}
\label{subsec:geometric-incidence-category}

\begin{definition}[Geometric incidence-completion category]
\label{def:geometric-incidence-category}
Let
\[
\mathbf{IncCov}(X,\Sigma)
\]
be the category whose objects are admissible incidence completions
\[
\overline\pi:
\overline Y\longrightarrow X
\]
of finite-sheeted coverings
\[
\pi^\circ:
Y^\circ\longrightarrow X^\circ.
\]

A morphism
\[
H:
(\overline Y,\overline\pi)
\longrightarrow
(\overline Z,\overline q)
\]
is a continuous map satisfying
\[
\overline q\circ H=\overline\pi.
\]
\end{definition}

Since
\[
(\overline\pi)^{-1}(X^\circ)=Y^\circ
\qquad\text{and}\qquad
(\overline q)^{-1}(X^\circ)=Z^\circ,
\]
every morphism restricts to a map
\[
H^\circ:
Y^\circ\longrightarrow Z^\circ
\]
over $X^\circ$.

Likewise, for each marked point $z_i$ it restricts to a continuous map
\[
\phi_i:
F_i\longrightarrow F_i',
\]
where
\[
F_i=\overline\pi^{-1}(z_i),
\qquad
F_i'=\overline q^{-1}(z_i).
\]

\begin{proposition}[Geometric morphism criterion]
\label{prop:geometric-morphism-criterion}
Let
\[
\overline Y
\qquad\text{and}\qquad
\overline Z
\]
be admissible incidence completions of finite coverings
\[
\pi^\circ:Y^\circ\to X^\circ
\]
and
\[
q^\circ:Z^\circ\to X^\circ.
\]

A map
\[
H:\overline Y\to\overline Z
\]
over $X$ is equivalent to the following data:

\begin{enumerate}[label=(\roman*)]
\item a map
\[
H^\circ:Y^\circ\to Z^\circ
\]
over $X^\circ$;

\item continuous maps
\[
\phi_i:F_i\to F_i';
\]

\item for every peripheral end $c$ of $Y^\circ$ above $z_i$,
\[
\phi_i(A_{i,c})
\subseteq
A_{i,H^\circ_*c}' ,
\]
where
\[
H^\circ_*c
\]
is the peripheral component of $Z^\circ$ containing
$H^\circ(E_{i,c})$.
\end{enumerate}
\end{proposition}

\begin{proof}
Suppose first that
\[
H:\overline Y\to\overline Z
\]
is continuous over $X$.  Its restrictions $H^\circ$ and $\phi_i$ are
continuous and lie over the appropriate parts of the base.

Since $E_{i,c}$ is connected, its image under $H^\circ$ lies in one
peripheral component of the target; denote that component by
\[
H^\circ_*c.
\]

Let
\[
a\in A_{i,c}.
\]
By exact endpoint realization, the $c$-tail converges to $a$.
Continuity of $H$ implies that the image tail converges to
\[
\phi_i(a).
\]
The image tail lies in the peripheral component $H^\circ_*c$, so exact
endpoint realization in the target gives
\[
\phi_i(a)\in A_{i,H^\circ_*c}'.
\]
Hence
\[
\phi_i(A_{i,c})
\subseteq
A_{i,H^\circ_*c}'.
\]

Conversely, suppose the stated data are given.  Define $H$ by
\[
H|_{Y^\circ}=H^\circ
\]
and
\[
H|_{F_i}=\phi_i.
\]

Continuity on $Y^\circ$ is already known.  Let
\[
a\in F_i
\]
and consider a basic neighbourhood
\[
B_i'(V',\mathbf r')
\]
of $\phi_i(a)$ in the target.

Set
\[
V:=\phi_i^{-1}(V').
\]
If
\[
c\in C_i(V),
\]
then some
\[
x\in A_{i,c}\cap V
\]
satisfies
\[
\phi_i(x)\in V'.
\]
By the incidence condition,
\[
\phi_i(x)
\in
A_{i,H^\circ_*c}'.
\]
Therefore
\[
H^\circ_*c
\in
C_i'(V').
\]

Because $H^\circ$ lies over $X^\circ$, sufficiently short $c$-tails are
mapped into sufficiently short $H^\circ_*c$-tails.  Since only finitely
many peripheral components occur, their radii can be chosen
simultaneously.  Hence a basic neighbourhood of $a$ is mapped into
$B_i'(V',\mathbf r')$.

Thus $H$ is continuous at every exceptional point.
\end{proof}

\subsection{The extraction functor}
\label{subsec:extraction-functor}

Let
\[
\overline\pi:\overline Y\to X
\]
be an object of
\[
\mathbf{IncCov}(X,\Sigma).
\]

Set
\[
A:=\overline\pi^{-1}(x_0).
\]
Since $x_0\in X^\circ$, this is the ordinary fibre of the regular
covering.  Path lifting gives it its canonical finite $G$-action.

For each $i$,
\[
H_i\backslash A
\]
is identified with the peripheral-end set $C_i$ by the standard orbit description of peripheral monodromy recorded above.

Set
\[
F_i:=\overline\pi^{-1}(z_i)
\]
with its subspace topology and define
\[
R_i
=
\left\{
(c,a):
a\in
\operatorname{Ends}_{\overline Y}(c)
\right\}.
\]
Exact endpoint realization identifies these endpoint sets with the
incidence sets used in constructing $\overline Y$.

Set
\[
\mathscr M(\overline Y)
\in
\mathbf{StrMon}(X,\Sigma,x_0)
\]
to be this structured datum.

If
\[
H:\overline Y\to\overline Z
\]
is a geometric morphism, its restriction to the regular fibre gives
\[
u:=H|_A:A\to B.
\]
Because $H^\circ$ lies over $X^\circ$, the classical monodromy
naturality of fibre maps implies that $u$ is $G$-equivariant.

The exceptional restrictions
\[
\phi_i:=H|_{F_i}
\]
satisfy the required incidence inclusions by
Proposition~\ref{prop:geometric-morphism-criterion}.

\begin{proposition}
\label{prop:extraction-functor}
The assignment above defines a functor
\[
\mathscr M:
\mathbf{IncCov}(X,\Sigma)
\longrightarrow
\mathbf{StrMon}(X,\Sigma,x_0).
\]
\end{proposition}

\subsection{Reconstruction from a finite \texorpdfstring{$G$}{G}-set}
\label{subsec:reconstruction-functor}

Conversely, let
\[
\mathfrak A
=
\left(
A;\{F_i,R_i\}
\right)
\]
be a structured monodromy datum.

The classical monodromy equivalence reconstructs a finite covering
\[
\pi_A^\circ:
Y_A^\circ\longrightarrow X^\circ
\]
whose fibre over $x_0$ is $A$ and whose monodromy action is the given
$G$-action.

One explicit construction uses the universal covering
\[
\widetilde X^\circ\longrightarrow X^\circ.
\]
With the usual compatible convention for the deck action of $G$, set
\[
Y_A^\circ
:=
\widetilde X^\circ\times_G A.
\]
Projection onto the first factor descends to
\[
\pi_A^\circ:Y_A^\circ\to X^\circ.
\]

For each $i$, its peripheral components are identified with
\[
H_i\backslash A.
\]
Use the given relation
\[
R_i
\subseteq
(H_i\backslash A)\times F_i
\]
to define endpoint sets
\[
A_{i,c}
=
\{a\in F_i:(c,a)\in R_i\}.
\]

Applying the incidence-completion construction of
Section~\ref{sec:incidence-completions} gives
\[
\overline Y_{\mathfrak A}
=
Y_A^\circ
\sqcup
\bigsqcup_iF_i
\]
and a continuous completed map
\[
\overline\pi_{\mathfrak A}:
\overline Y_{\mathfrak A}\to X.
\]

Now let
\[
(u,\phi_i):
\mathfrak A\to\mathfrak B
\]
be a structured morphism.

The $G$-map
\[
u:A\to B
\]
induces, by the classical monodromy equivalence, a unique map of
coverings
\[
u^\circ:
Y_A^\circ\to Y_B^\circ
\]
over $X^\circ$.

Its action on peripheral components is exactly the orbit map
\[
\overline u_i:
H_i\backslash A
\to
H_i\backslash B.
\]

Since
\[
\phi_i(A_{i,c})
\subseteq
A_{i,\overline u_i(c)}',
\]
Proposition~\ref{prop:geometric-morphism-criterion} extends
\[
u^\circ
\]
and the $\phi_i$ uniquely to a morphism
\[
\overline u:
\overline Y_{\mathfrak A}
\to
\overline Y_{\mathfrak B}
\]
over $X$.

Therefore:

\begin{proposition}
\label{prop:reconstruction-functor}
The construction above defines a functor
\[
\mathscr R:
\mathbf{StrMon}(X,\Sigma,x_0)
\longrightarrow
\mathbf{IncCov}(X,\Sigma).
\]
\end{proposition}

\subsection{Equivalence theorem}
\label{subsec:structured-equivalence}

\begin{theorem}[Structured monodromy classification]
\label{thm:structured-monodromy-equivalence}
There is an equivalence of categories
\[
\mathbf{IncCov}(X,\Sigma)
\simeq
\mathbf{StrMon}(X,\Sigma,x_0).
\]

Under this equivalence, an incidence completion of a finite-sheeted regular covering is
determined by its finite $G$-set together with the exceptional fibres and their peripheral
endpoint incidence.
\end{theorem}

\begin{proof}
The classical monodromy equivalence identifies the category of finite
coverings of $X^\circ$ with the category of finite $G$-sets.  Therefore,
for every finite $G$-set $A$,
\[
\mathscr M
\bigl(
\mathscr R(\mathfrak A)
\bigr)
\]
recovers $A$, up to the canonical isomorphism supplied by that
equivalence.

The exceptional fibre $F_i$ is inserted unchanged by
$\mathscr R$.  Exact endpoint realization then recovers from the
completed topology precisely the prescribed endpoint sets
\[
A_{i,c}.
\]
Hence
\[
\mathscr M\circ\mathscr R
\cong
\operatorname{id}_{\mathbf{StrMon}}.
\]

Conversely, begin with a geometric incidence completion
\[
\overline Y.
\]
Extraction produces its regular monodromy $G$-set
\[
A=\overline\pi^{-1}(x_0)
\]
and its intrinsic endpoint relations.

Reconstruction from $A$ gives a regular covering canonically isomorphic,
under classical monodromy theory, to the original regular covering
$Y^\circ\to X^\circ$.  Under this regular isomorphism the peripheral
components correspond.  Since the exceptional fibres and all endpoint
sets are identical, Proposition~\ref{prop:geometric-morphism-criterion} extends the
regular isomorphism, together with the identity maps on the exceptional fibres, to an
isomorphism
\[
\mathscr R\bigl(\mathscr M(\overline Y)\bigr)
\cong
\overline Y
\]
over $X$.

Naturality on morphisms follows from
Proposition~\ref{prop:geometric-morphism-criterion} and the full
faithfulness of the ordinary monodromy equivalence.

Thus
\[
\mathscr R\circ\mathscr M
\cong
\operatorname{id}_{\mathbf{IncCov}}
\]
and
\[
\mathscr M\circ\mathscr R
\cong
\operatorname{id}_{\mathbf{StrMon}},
\]
proving the equivalence.
\end{proof}

\subsection{Isomorphism criterion}
\label{subsec:structured-isomorphism}

\begin{corollary}[Structured isomorphism criterion]
\label{cor:structured-isomorphism}
Two incidence completions are isomorphic over $X$ if and only if there
exist:

\begin{enumerate}[label=(\roman*)]
\item a $G$-set isomorphism
\[
u:A\xrightarrow{\ \cong\ }B;
\]

\item homeomorphisms
\[
\phi_i:F_i\xrightarrow{\ \cong\ }F_i';
\]

\item for every peripheral orbit
\[
c\in H_i\backslash A,
\]
the equality
\[
\phi_i(A_{i,c})
=
A_{i,\overline u_i(c)}'.
\]
\end{enumerate}
\end{corollary}

\section{Automorphisms and the singular symmetry kernel}
\label{sec:singular-automorphisms}

Fix an incidence completion $\overline\pi:\overline Y\to X$ with regular
restriction $\pi^\circ:Y^\circ\to X^\circ$, and let
$\mathfrak A=(A;\{F_i,R_i\}_{i=1}^n)$ be its structured monodromy datum.  With
$G=\pi_1(X^\circ,x_0)$, the automorphism group of the regular covering is
canonically $\operatorname{Aut}_G(A)$.

The completed automorphism group is
\[
\operatorname{Aut}_X(\overline Y)
:=
\left\{
H\in\operatorname{Homeo}(\overline Y):
\overline\pi\circ H=\overline\pi
\right\}.
\]

\subsection{Restriction to the regular covering}
\label{subsec:restriction-regular}

Every automorphism of the completed space preserves the regular locus
because
\[
Y^\circ
=
\overline\pi^{-1}(X^\circ).
\]
Hence restriction gives a homomorphism
\[
\operatorname{res}:
\operatorname{Aut}_X(\overline Y)
\longrightarrow
\operatorname{Aut}_{X^\circ}(Y^\circ).
\]

Via classical monodromy,
\[
\operatorname{Aut}_{X^\circ}(Y^\circ)
\cong
\operatorname{Aut}_G(A).
\]

Not every regular automorphism need extend across the exceptional fibres.
The incidence relation may break part of the ordinary covering symmetry.

\begin{definition}[Admissible regular automorphisms]
\label{def:admissible-regular-aut}
Define
\[
\operatorname{Aut}_G(A)_{\mathrm{adm}}
\]
to be the subgroup of all
\[
u\in\operatorname{Aut}_G(A)
\]
for which there exist homeomorphisms
\[
\phi_i:F_i\longrightarrow F_i
\]
such that
\[
\phi_i(A_{i,c})
=
A_{i,\overline u_i(c)}
\]
for every $i$ and every
\[
c\in H_i\backslash A,
\]
where
\[
\overline u_i:
H_i\backslash A
\longrightarrow
H_i\backslash A
\]
is the induced permutation of peripheral orbits.
\end{definition}

\begin{proposition}
\label{prop:admissible-subgroup}
The set
\[
\operatorname{Aut}_G(A)_{\mathrm{adm}}
\]
is a subgroup of
\[
\operatorname{Aut}_G(A).
\]
Moreover,
\[
\operatorname{im}(\operatorname{res})
=
\operatorname{Aut}_G(A)_{\mathrm{adm}}.
\]
\end{proposition}

\begin{proof}
The identity is admissible by taking every $\phi_i$ to be the identity.

Suppose
\[
u,v\in\operatorname{Aut}_G(A)_{\mathrm{adm}}
\]
are realized by fibre homeomorphisms
\[
\phi_i
\qquad\text{and}\qquad
\psi_i.
\]
Then
\[
\phi_i(A_{i,c})
=
A_{i,\overline u_i(c)}
\]
and
\[
\psi_i(A_{i,d})
=
A_{i,\overline v_i(d)}.
\]
Therefore
\[
(\psi_i\circ\phi_i)(A_{i,c})
=
A_{i,
\overline v_i(\overline u_i(c))}
=
A_{i,\overline{(vu)}_i(c)}.
\]
Thus $vu$ is admissible.

Likewise,
\[
\phi_i^{-1}
\left(
A_{i,\overline u_i(c)}
\right)
=
A_{i,c},
\]
so $u^{-1}$ is admissible.

The equality with the image of restriction follows from the structured
isomorphism criterion:
an automorphism of the completed space restricts to a $G$-set
automorphism satisfying precisely the displayed compatibility, while any
such compatible pair extends to an automorphism of the incidence
completion.
\end{proof}

The incidence datum may reduce the ordinary deck symmetry:
\[
\operatorname{Aut}_G(A)_{\mathrm{adm}}
\leq
\operatorname{Aut}_G(A).
\]

In the rigid symmetric case equality may hold; for asymmetric endpoint
incidence it can be strict.

\subsection{The singular automorphism kernel}
\label{subsec:singular-kernel}

\begin{definition}[Singular automorphism kernel]
\label{def:singular-kernel}
The \emph{singular automorphism kernel} of the incidence completion is
\[
K_{\mathrm{sing}}(\overline Y)
:=
\ker
\left(
\operatorname{Aut}_X(\overline Y)
\longrightarrow
\operatorname{Aut}_{X^\circ}(Y^\circ)
\right).
\]
\end{definition}

An element of $K_{\mathrm{sing}}$ fixes the entire regular locus as a map
over $X^\circ$ and acts only on the exceptional fibres.

For each $i$, define
\[
\operatorname{Aut}(F_i;\mathcal A_i)
:=
\left\{
\phi\in\operatorname{Homeo}(F_i):
\phi(A_{i,c})=A_{i,c}
\text{ for all }c\in C_i
\right\}.
\]

\begin{theorem}[Singular-kernel theorem]
\label{thm:singular-kernel}
There is a natural isomorphism
\[
K_{\mathrm{sing}}(\overline Y)
\cong
\prod_{i=1}^n
\operatorname{Aut}(F_i;\mathcal A_i).
\]
\end{theorem}

\begin{proof}
Let
\[
H\in K_{\mathrm{sing}}(\overline Y).
\]
By definition,
\[
H|_{Y^\circ}
=
\operatorname{id}_{Y^\circ}.
\]
For each $i$, its restriction
\[
\phi_i:=H|_{F_i}
\]
is a homeomorphism of $F_i$.

Since the induced regular automorphism is the identity, every peripheral
end $c$ is fixed.  The structured isomorphism criterion therefore gives
\[
\phi_i(A_{i,c})=A_{i,c}
\]
for every $c$.  Hence
\[
\phi_i\in\operatorname{Aut}(F_i;\mathcal A_i).
\]

Conversely, given a tuple
\[
(\phi_i)_i
\in
\prod_i\operatorname{Aut}(F_i;\mathcal A_i),
\]
extend it by the identity on $Y^\circ$.  The endpoint equalities
\[
\phi_i(A_{i,c})=A_{i,c}
\]
satisfy Corollary~\ref{cor:structured-isomorphism}, so the resulting set map is an
automorphism of $\overline Y$ over $X$.

The two constructions are inverse and respect composition.
\end{proof}

\subsection{The automorphism exact sequence}
\label{subsec:automorphism-exact-sequence}

\begin{theorem}[Automorphism exact sequence]
\label{thm:completed-aut-exact}
There is a natural short exact sequence
\[
1
\longrightarrow
K_{\mathrm{sing}}(\overline Y)
\longrightarrow
\operatorname{Aut}_X(\overline Y)
\overset{\operatorname{res}}{\longrightarrow}
\operatorname{Aut}_G(A)_{\mathrm{adm}}
\longrightarrow
1.
\]
Equivalently,
\[
1
\longrightarrow
\prod_i
\operatorname{Aut}(F_i;\mathcal A_i)
\longrightarrow
\operatorname{Aut}_X(\overline Y)
\longrightarrow
\operatorname{Aut}_G(A)_{\mathrm{adm}}
\longrightarrow
1.
\]
\end{theorem}

\begin{proof}
The kernel is
Theorem~\ref{thm:singular-kernel}, while
Proposition~\ref{prop:admissible-subgroup} identifies the image with
$\operatorname{Aut}_G(A)_{\mathrm{adm}}$.
\end{proof}

No splitting is asserted; a semidirect-product description requires additional
hypotheses ensuring compatible exceptional lifts of admissible regular
automorphisms.

\subsection{Finite fibres and weighted-poset automorphisms}
\label{subsec:finite-kernel-posets}

Assume now that every exceptional fibre is finite and let
\[
(P_i,\mu_i,\mathcal I_i)
\]
be its finite incidence triple.

Define
\[
\operatorname{Aut}(P_i,\mu_i,\mathcal I_i)
\]
to be the group of all poset automorphisms
\[
\alpha:P_i\to P_i
\]
satisfying
\[
\mu_i(\alpha(p))=\mu_i(p)
\]
and
\[
\mathcal I_i(\alpha(p))
=
\mathcal I_i(p)
\]
for every $p\in P_i$.

An incidence-preserving homeomorphism of $F_i$ induces such a poset
automorphism.

The kernel of this induced action consists exactly of arbitrary
permutations inside the topological indistinguishability classes.

\begin{proposition}[Finite singular-kernel sequence]
\label{prop:finite-singular-kernel-sequence}
For every finite exceptional fibre there is a natural short exact
sequence
\[
1
\longrightarrow
\prod_{p\in P_i}
\Sym_{\mu_i(p)}
\longrightarrow
\operatorname{Aut}(F_i;\mathcal A_i)
\longrightarrow
\operatorname{Aut}(P_i,\mu_i,\mathcal I_i)
\longrightarrow
1.
\]
Consequently
\[
K_{\mathrm{sing}}
\]
is assembled from permutations inside indistinguishability classes and
symmetries of the weighted specialization posets preserving incidence.
\end{proposition}

\begin{proof}
A homeomorphism of a finite space preserves the specialization preorder
and therefore induces an automorphism of the Kolmogorov quotient.
Because it is a bijection, it preserves the sizes of the
indistinguishability classes.

Preservation of every endpoint set is equivalent to preservation of the
incidence labels
\[
\mathcal I_i(p).
\]
Thus the induced automorphism lies in
\[
\operatorname{Aut}(P_i,\mu_i,\mathcal I_i).
\]

If the induced poset automorphism is the identity, the homeomorphism may
permute points only within each indistinguishability class.  Every such
permutation is a homeomorphism and preserves all endpoint sets, because
the latter are saturated.  Hence the kernel is
\[
\prod_{p\in P_i}\Sym_{\mu_i(p)}.
\]

Conversely, given
\[
\alpha\in
\operatorname{Aut}(P_i,\mu_i,\mathcal I_i),
\]
weight preservation allows a bijection
\[
F_p\longrightarrow F_{\alpha(p)}
\]
for each specialization class.  Combining these bijections gives a
homeomorphism of $F_i$ inducing $\alpha$.  Since $\alpha$ preserves
$\mathcal I_i$, the homeomorphism preserves every endpoint set.
Therefore the final map is surjective.
\end{proof}

\subsection{Quotient-type intermediate completions}
\label{subsec:quotient-intermediate}

Let
\[
\overline\pi:
\overline Y\longrightarrow X
\]
be the incidence completion corresponding to $\mathfrak A$.

\begin{definition}[Quotient-type intermediate completion]
\label{def:quotient-intermediate}
A \emph{quotient-type intermediate completion} of $\overline Y$ is a
factorization
\[
\overline Y
\overset{Q}{\longrightarrow}
\overline Z
\overset{\overline q}{\longrightarrow}
X
\]
such that:
\begin{enumerate}[label=(\roman*)]
\item
\[
\overline q\circ Q=\overline\pi;
\]
\item $Q$ is surjective;
\item its restriction
\[
Q^\circ:Y^\circ\to Z^\circ
\]
is a quotient morphism of finite coverings over $X^\circ$;
\item for every $i$, the induced map
\[
Q_i:F_i\to F_i'
\]
on exceptional fibres is a surjective quotient map.
\end{enumerate}
Two such factorizations are equivalent if their intermediate objects are
isomorphic over $X$ compatibly with the quotient maps from
$\overline Y$.
\end{definition}

The restriction
\[
Q^\circ:Y^\circ\to Z^\circ
\]
corresponds under ordinary monodromy theory to a surjective $G$-map
\[
q:A\longrightarrow B.
\]

For every puncture $z_i$, this induces a surjection on peripheral orbit
sets
\[
\overline q_i:
H_i\backslash A
\longrightarrow
H_i\backslash B.
\]

Ordinary intermediate covering theory therefore supplies the quotient
$A\twoheadrightarrow B$.

\subsection{Structured quotient data}
\label{subsec:structured-quotient-data}

\begin{definition}[Structured quotient datum]
\label{def:structured-quotient-datum}
A \emph{structured quotient datum} of $\mathfrak A$ consists of:
\begin{enumerate}[label=(\roman*)]
\item a surjective $G$-map
\[
q:A\twoheadrightarrow B;
\]
\item for every $i$, a surjective quotient map of topological spaces
\[
\psi_i:F_i\twoheadrightarrow F_i';
\]
\item for every
\[
d\in H_i\backslash B,
\]
a nonempty closed subset
\[
B_{i,d}\subseteq F_i'
\]
such that
\[
F_i'
=
\bigcup_{d\in H_i\backslash B}B_{i,d}
\]
and
\[
\psi_i(A_{i,c})
\subseteq
B_{i,\overline q_i(c)}
\]
for every
\[
c\in H_i\backslash A.
\]
\end{enumerate}
\end{definition}

The last condition says that every endpoint of a source peripheral end
maps to an allowed endpoint of the corresponding quotient peripheral
end.

The target datum
\[
\mathfrak B
=
\left(
B;\{F_i',B_{i,d}\}
\right)
\]
is itself a structured monodromy datum, and
\[
(q,\psi_i):
\mathfrak A\to\mathfrak B
\]
is a surjective structured morphism.

\begin{proposition}
\label{prop:structured-quotient-geometry}
Structured quotient data of $\mathfrak A$ are equivalent, under the
classification of
Theorem~\ref{thm:structured-monodromy-equivalence}, to quotient-type
intermediate completions of $\overline Y$.
\end{proposition}

\begin{proof}
Given structured quotient data, the structured morphism
\[
(q,\psi_i):
\mathfrak A\to\mathfrak B
\]
corresponds under
Theorem~\ref{thm:structured-monodromy-equivalence}
to a morphism
\[
Q:\overline Y\to\overline Z
\]
over $X$.

The maps $q$ and $\psi_i$ are surjective, so $Q$ is surjective on the
regular locus and on every exceptional fibre, hence on the total space.
Its regular restriction is the quotient morphism corresponding to the
surjective $G$-map $q$.

Conversely, a quotient-type intermediate completion yields a surjective
map of regular fibres
\[
q:A\twoheadrightarrow B
\]
and surjective exceptional-fibre maps
\[
\psi_i:F_i\twoheadrightarrow F_i'.
\]
The geometric morphism criterion gives
\[
\psi_i(A_{i,c})
\subseteq
B_{i,\overline q_i(c)}.
\]
Thus the extracted data form a structured quotient datum.

The two constructions are inverse up to the stated equivalence of
intermediate factorizations.
\end{proof}

\subsection{Canonical minimal quotient incidence}
\label{subsec:minimal-quotient-incidence}

Fix only the quotient maps
\[
q:A\twoheadrightarrow B
\]
and
\[
\psi_i:F_i\twoheadrightarrow F_i'.
\]

For
\[
d\in H_i\backslash B,
\]
define the raw endpoint image
\[
B_{i,d}^{\mathrm{raw}}
:=
\bigcup_{\substack{c\in H_i\backslash A\\
\overline q_i(c)=d}}
\psi_i(A_{i,c}).
\]

This set is nonempty because $\overline q_i$ is surjective and every
source endpoint set is nonempty.

Define
\[
B_{i,d}^{\min}
:=
\overline{
B_{i,d}^{\mathrm{raw}}
}^{\,F_i'}.
\]

\begin{proposition}[Minimal quotient incidence]
\label{prop:minimal-quotient-incidence}
The family
\[
\{B_{i,d}^{\min}\}_{d\in H_i\backslash B}
\]
is an admissible closed incidence family on $F_i'$.

Moreover, if
\[
\{B_{i,d}\}
\]
is any closed incidence family making
\[
(q,\psi_i)
\]
a structured morphism, then
\[
B_{i,d}^{\min}
\subseteq
B_{i,d}
\]
for every $i,d$.
\end{proposition}

\begin{proof}
Every $B_{i,d}^{\min}$ is closed and nonempty by construction.

To prove endpoint-covering, let
\[
b\in F_i'.
\]
Choose
\[
a\in F_i
\]
with
\[
\psi_i(a)=b.
\]
Since the source datum is endpoint-covering, there exists
\[
c\in H_i\backslash A
\]
such that
\[
a\in A_{i,c}.
\]
Then
\[
b=\psi_i(a)
\in
B_{i,\overline q_i(c)}^{\mathrm{raw}}
\subseteq
B_{i,\overline q_i(c)}^{\min}.
\]
Hence the minimal endpoint sets cover $F_i'$.

Now let $\{B_{i,d}\}$ be any compatible closed target incidence family.
For every source orbit $c$,
\[
\psi_i(A_{i,c})
\subseteq
B_{i,\overline q_i(c)}.
\]
Therefore
\[
B_{i,d}^{\mathrm{raw}}
\subseteq
B_{i,d}.
\]
Since $B_{i,d}$ is closed,
\[
B_{i,d}^{\min}
=
\overline{B_{i,d}^{\mathrm{raw}}}
\subseteq
B_{i,d}.
\]
\end{proof}

The quotient maps therefore force at least the incidence
\[
B_{i,d}^{\min},
\]
but an intermediate completion may possess larger endpoint sets.

\subsection{Structured quotient correspondence}
\label{subsec:structured-quotient-correspondence}

\begin{theorem}[Structured quotient correspondence]
\label{thm:structured-quotient-correspondence}
Let $\overline Y\to X$ be an incidence completion with structured monodromy datum
$\mathfrak A$.  Equivalence classes of quotient-type intermediate completions
\[
\overline Y\twoheadrightarrow\overline Z\longrightarrow X
\]
are naturally in bijection with isomorphism classes of structured quotient data of
$\mathfrak A$.  The regular part is a finite $G$-set quotient $A\twoheadrightarrow B$;
the exceptional fibres are quotient spaces $F_i\twoheadrightarrow F_i'$; and each target
endpoint set is a closed enlargement of the canonical minimal image incidence.
\end{theorem}

\begin{proof}
This is Proposition~\ref{prop:structured-quotient-geometry} together with
Proposition~\ref{prop:minimal-quotient-incidence}.
\end{proof}

\begin{corollary}[Rigid sector]
Assume that each source exceptional fibre is discrete and that its endpoint relation is
the graph of a bijection $H_i\backslash A\xrightarrow{\sim}F_i$.  Restrict quotient
targets to the same rigid form: for a quotient $A\twoheadrightarrow B$, require
$F_i'=H_i\backslash B$ with the discrete topology and $B_{i,d}=\{d\}$ for every
$d\in H_i\backslash B$.  In this sector the exceptional data are determined by the
finite $G$-set quotient.  If $A\cong G/K$ is connected, rigid intermediate completions
reduce to the classical correspondence with intermediate subgroups $K\leq L\leq G$.
\end{corollary}

\section{Singular transport through incidence defects}
\label{sec:local-singular-transport}

Regular paths in $X^\circ$ act by the usual covering-transport bijections.
For an isolated crossing of $z_i\in\Sigma$, continuation from an incoming
peripheral germ to an outgoing germ is determined by their common
exceptional endpoints.

\subsection{Incoming and outgoing peripheral germs}
\label{subsec:incoming-outgoing-germs}

Fix $z_i\in\Sigma$ and suppress the index $i$.  Let $C$ be the peripheral-end
set and $F$ the exceptional fibre, with endpoint sets $A_c\subseteq F$ for $c\in C$.

Consider a path $\gamma:[-1,1]\to X$ with $\gamma(0)=z_i$ and
$\gamma(t)\in X^\circ$ for $t\neq0$.

Suppose a regular lift of the negative half approaches the defect through
a peripheral end
\[
b\in C,
\]
while a regular lift of the positive half leaves through
\[
c\in C.
\]

The accumulated-end set of the resulting punctured lift at $0$ is
\[
\{b,c\},
\]
with the understanding that it reduces to $\{b\}$ when $b=c$.

\begin{definition}[Local transport set]
\label{def:local-transport-set}
For
\[
b,c\in C,
\]
define the \emph{local transport set}
\[
T(b,c)
:=
A_b\cap A_c.
\]
\end{definition}

A crossing from $b$ to $c$ is realizable precisely when
\[
T(b,c)\neq\varnothing.
\]

\begin{proposition}[Two-sided continuation criterion]
\label{prop:two-sided-continuation}
Let $\gamma$ and its regular lift be as above.  Then the possible values
of a continuous completed lift at the defect time $0$ are exactly
\[
T(b,c)=A_b\cap A_c.
\]
\end{proposition}

\begin{proof}
The accumulated-end set at $0$ is
\[
\Lambda(0)=\{b,c\}.
\]
The statement is therefore the case
\[
J=\{b,c\}
\]
of Theorem~\ref{thm:incidence-lift}.
\end{proof}

For finite discrete $F$, the crossing is obstructed, unique, or multiply realized according as $|T(b,c)|$ is $0$, $1$, or greater than $1$.

\subsection{Multiplicity-valued transport}
\label{subsec:multiplicity-transport}

Assume for the remainder of this subsection that $F$ is finite and
discrete.

\begin{definition}[Transport multiplicity matrix]
\label{def:transport-matrix}
Define
\[
M=(M_{bc})_{b,c\in C}
\]
by
\[
M_{bc}
:=
|A_b\cap A_c|.
\]
\end{definition}

Then
\[
M_{bc}
\]
is the number of admissible exceptional values for passage from $b$ to $c$.

The matrix is symmetric:
\[
M_{bc}=M_{cb},
\]
and its diagonal entries are
\[
M_{cc}=|A_c|.
\]

Hence the diagonal records endpoint multiplicities of individual
peripheral ends, while the off-diagonal entries record pairwise
end-merging multiplicities.

Let
\[
B
\]
be the $|C|\times |F|$ incidence matrix
\[
B_{ca}
=
\begin{cases}
1,&a\in A_c,\\
0,&a\notin A_c.
\end{cases}
\]

\begin{proposition}[Gram factorization]
\label{prop:transport-Gram}
The transport multiplicity matrix satisfies
\[
M=BB^{\mathsf T}.
\]
\end{proposition}

\begin{proof}
The $(b,c)$ entry of
\[
BB^{\mathsf T}
\]
is
\[
\sum_{a\in F}B_{ba}B_{ca},
\]
which counts exactly the points
\[
a\in A_b\cap A_c.
\]
\end{proof}

Thus $M$ is positive semidefinite over $\mathbb R$.

However, the factorization does not imply that $M$ determines $B$.
As shown in Section~\ref{sec:lift-reconstruction}, pairwise intersection
counts do not determine higher endpoint incidence.

\subsection{Higher singular transport}
\label{subsec:higher-singular-transport}

Let
\[
J\subseteq C
\]
be a nonempty finite set of peripheral ends.

\begin{definition}[Higher transport set]
\label{def:higher-transport-set}
Define
\[
T(J)
:=
\bigcap_{c\in J}A_c.
\]
\end{definition}

\begin{proposition}[Higher continuation criterion]
\label{prop:higher-continuation}
Suppose a regular lift accumulates at one defect point through precisely
the peripheral ends in
\[
J\subseteq C.
\]
Then its possible singular values are exactly
\[
T(J).
\]
\end{proposition}

Thus the full local transport structure is the contravariant diagram
\[
T:
\mathcal P(C)^{\mathrm{op}}
\longrightarrow
\mathbf{Top},
\qquad
J\longmapsto T(J),
\]
with
\[
T(\{c\})=A_c.
\]

For
\[
J\subseteq K,
\]
there is a natural inclusion
\[
T(K)\hookrightarrow T(J).
\]

\begin{example}[Same pairwise transport, different higher transport]
\label{ex:same-pairwise-different-higher}
Let
\[
C=\{1,2,3\}.
\]

In the first incidence datum take three exceptional points
\[
a_{12},a_{13},a_{23}
\]
with
\[
A_1=\{a_{12},a_{13}\},
\]
\[
A_2=\{a_{12},a_{23}\},
\]
\[
A_3=\{a_{13},a_{23}\}.
\]
Then every pairwise transport set has cardinality one, but
\[
T(\{1,2,3\})=\varnothing.
\]

In the second datum, take a single exceptional point $a$ and put
\[
A_1=A_2=A_3=\{a\}.
\]
Again every pairwise transport set has cardinality one, but now
\[
T(\{1,2,3\})=\{a\}.
\]

Hence the pairwise multiplicity matrix is identical in the two examples,
while the three-arm transport differs.
\end{example}

\subsection{Decorated finite-crossing category}
\label{subsec:finite-crossing-category}

For each $z_i\in\Sigma$, choose a sufficiently small coordinate disc $U_i$,
a reference point $u_i\in U_i^\times$, and a radial arc
\[
\rho_i:[0,1]\longrightarrow U_i,
\qquad
\rho_i(0)=u_i,\quad
\rho_i(1)=z_i,\quad
\rho_i([0,1))\subseteq U_i^\times.
\]
Let $A_i^{\mathrm{reg}}=(\pi^\circ)^{-1}(u_i)$ and let
$\kappa_i:A_i^{\mathrm{reg}}\to C_i$ send a sheet to its peripheral
component germ.

\begin{definition}[Decorated finite-crossing category]
\label{def:finite-crossing-category}
Let $\mathcal P_{\mathrm{fc}}(X,\Sigma;\mathbf u)$ be the category generated by
the fundamental groupoid $\Pi_1(X^\circ)$ together with, for each $i$, one
additional morphism
\[
\delta_i:u_i\longrightarrow u_i,
\]
subject to no relations involving the $\delta_i$ beyond the category axioms.
With composition read right-to-left, every morphism has an alternating
representative
\[
[\eta_m]\,\delta_{i_m}\,[\eta_{m-1}]\cdots
\delta_{i_1}\,[\eta_0],
\]
where adjacent sources and targets agree.
\end{definition}

Define the fixed crossing path
\[
\chi_i(t)
=
\begin{cases}
\rho_i(2t),&0\leq t\leq \frac12,\\[1mm]
\rho_i(2-2t),&\frac12\leq t\leq1.
\end{cases}
\]
It meets $\Sigma$ only at $t=\frac12$ and represents the geometric crossing
assigned to $\delta_i$.

\begin{definition}[Standard geometric representative]
\label{def:standard-geometric-representative}
For an alternating word
\[
w=[\eta_m]\,\delta_{i_m}\,[\eta_{m-1}]\cdots
\delta_{i_1}\,[\eta_0],
\]
a \emph{standard geometric representative} is obtained by choosing
endpoint-fixed representatives of the regular classes $[\eta_j]$, inserting
$\chi_{i_j}$ for each $\delta_{i_j}$, and reparametrizing the concatenation
so that the singular crossing times are distinct.
\end{definition}

\begin{proposition}[Change of local crossing data]
\label{prop:change-crossing-data}
Let $(u_i',\rho_i',\delta_i')$ arise from another local choice in the same
punctured coordinate neighbourhood, and let
$\kappa_i':(\pi^\circ)^{-1}(u_i')\to C_i$ be the corresponding germ map.
Choose a regular path class $h_i:u_i\to u_i'$ represented in $U_i^\times$.
Ordinary covering transport $P_{h_i}$ then satisfies
\[
\kappa_i'\circ P_{h_i}=\kappa_i.
\]
The identity on $\Pi_1(X^\circ)$ extends to an isomorphism of generated
categories by
\[
\delta_i\longmapsto
[h_i^{-1}]\,\delta_i'\,[h_i].
\]
Under ordinary covering transport along $h_i$, the corresponding
relation- and multiplicity-valued transport functors are naturally
identified.
\end{proposition}

\begin{proof}
The displayed word is an endomorphism of $u_i$, so the universal property
of the freely generated category gives a functor fixing
$\Pi_1(X^\circ)$.  The inverse sends
\[
\delta_i'\longmapsto
[h_i]\,\delta_i\,[h_i^{-1}].
\]
Ordinary path lifting along $h_i$ identifies the local regular fibres and
conjugates the crossing relations and multiplicity matrices by the induced
bijection.
\end{proof}

Woolf's fundamental category treats homotopically stratified spaces and
classifies stratified \mbox{\'etale} covers and, dually, stratified branched
covers under the hypotheses of that theory \citep{Woolf2009}.  The line with
two origins appears there as a stratified \mbox{\'etale} cover.  The category
$\mathcal P_{\mathrm{fc}}(X,\Sigma;\mathbf u)$ is a narrower generated category
for isolated crossings and is used only to encode continuation determined by
the endpoint-incidence datum; no identification with Woolf's fundamental,
exit-, or entrance-path category is asserted.

\subsection{Sheet-level singular crossing relation}
\label{subsec:sheet-level-crossings}

Define
\[
\mathcal S_i\subseteq A_i^{\mathrm{reg}}\times A_i^{\mathrm{reg}}
\]
by
\[
(y_-,y_+)\in\mathcal S_i
\quad\Longleftrightarrow\quad
A_{i,\kappa_i(y_-)}\cap A_{i,\kappa_i(y_+)}\neq\varnothing.
\]
For finite discrete $F_i$, define the multiplicity matrix
\[
(\mathcal M_i)_{y_+,y_-}
:=\left|A_{i,\kappa_i(y_-)}\cap A_{i,\kappa_i(y_+)}\right|.
\]
Its support is $\mathcal S_i$.  At the end level this is the matrix
$B_iB_i^{\mathsf T}$, where $B_i$ is the incidence matrix from $C_i$ to $F_i$.

\begin{proposition}[Geometric meaning of the singular generator]
\label{prop:crossing-relation-geometric}
A pair $(y_-,y_+)\in A_i^{\mathrm{reg}}\times A_i^{\mathrm{reg}}$ lies in
$\mathcal S_i$ if and only if the two standard radial regular lifts determined by
$y_-$ and $y_+$ extend through $z_i$ to a continuous lift of the crossing represented by
$\delta_i$.
\end{proposition}

\begin{proof}
The incoming and outgoing pieces accumulate through the germs
$\kappa_i(y_-)$ and $\kappa_i(y_+)$.  The claim is the two-germ case of
Theorem~\ref{thm:incidence-lift}.
\end{proof}

\subsection{Relation-valued transport}
\label{subsec:crossing-word-transport}

Let $A_x=(\pi^\circ)^{-1}(x)$ for $x\in X^\circ$, and let $\mathbf{FinRel}$ denote the category of finite sets and relations, with relational composition.  For a regular path class $[\eta]:x\to y$, let $P_\eta:A_x\to A_y$ be ordinary
covering transport and let $\Gamma(P_\eta)$ be its graph.  Define
\[
\mathscr T([\eta])=\Gamma(P_\eta),
\qquad
\mathscr T(\delta_i)=\mathcal S_i,
\]
and extend by relational composition to every morphism of
$\mathcal P_{\mathrm{fc}}(X,\Sigma;\mathbf u)$.

\begin{theorem}[Geometric realization of incidence transport]
\label{thm:transport-realization}
Let $w:x\to y$ be a morphism of
$\mathcal P_{\mathrm{fc}}(X,\Sigma;\mathbf u)$ and let $a\in A_x$, $b\in A_y$.
Then $(a,b)\in\mathscr T(w)$ if and only if some, equivalently every,
standard geometric representative of $w$ admits a continuous lift to
$\overline Y$ beginning at $a$ and ending at $b$.
\end{theorem}

\begin{proof}
Changing an endpoint-fixed representative of a regular path class does not change
its covering transport, by the covering homotopy property.  Hence lift existence is
independent of the choices made in Definition~\ref{def:standard-geometric-representative}.
Write $w$ as an alternating word in regular path classes and singular generators.  Along
a regular factor, the lift is the unique ordinary covering lift.  At a singular factor
$\delta_i$, Proposition~\ref{prop:crossing-relation-geometric} identifies the admissible
incoming and outgoing sheet states.  Relational composition therefore records exactly the
choices of intermediate sheet states for which every factor lifts.  The singular times are
distinct in the standard representative, so the corresponding exceptional values can be
chosen independently.  The converse follows by restricting any completed lift to the
regular and singular factors.
\end{proof}

\begin{theorem}[Incidence transport functor]
\label{thm:incidence-transport-functor}
The assignment above defines a functor
\[
\mathscr T:
\mathcal P_{\mathrm{fc}}(X,\Sigma;\mathbf u)
\longrightarrow\mathbf{FinRel}.
\]
Its restriction to $\Pi_1(X^\circ)$ is the ordinary finite-cover transport functor.
\end{theorem}

\begin{proof}
The category $\mathcal P_{\mathrm{fc}}$ is generated by $\Pi_1(X^\circ)$ and the
$\delta_i$.  Ordinary transport respects composition in $\Pi_1(X^\circ)$, and
relational composition is associative.  The universal property of the generated category
therefore gives the functor.
\end{proof}

\subsection{Multiplicity-valued sheet transport}
\label{subsec:word-multiplicity}

Let $\mathbf{Mat}_{\mathbb N}$ denote the category of finite sets and matrices with entries in $\mathbb N$, with composition by matrix multiplication.

For an alternating word
\[
w=[\eta_m]\,\delta_{i_m}\,[\eta_{m-1}]\cdots\delta_{i_1}\,[\eta_0],
\]
let $P_{\eta_j}$ also denote the corresponding permutation matrix and define
\[
\mathscr T_{\mathbb N}(w)=P_{\eta_m}\mathcal M_{i_m}P_{\eta_{m-1}}\cdots\mathcal M_{i_1}P_{\eta_0}.
\]

This is a matrix
\[
A_x\longrightarrow A_y
\]
over $\mathbb N$.

\begin{theorem}[Multiplicity transport theorem]
\label{thm:multiplicity-transport}
Assume all exceptional fibres are finite and discrete.  Let
\[
w:x\to y
\]
be a morphism of $\mathcal P_{\mathrm{fc}}(X,\Sigma;\mathbf u)$, and let
\[
a\in A_x,
\qquad
b\in A_y.
\]

Then
\[
\mathscr T_{\mathbb N}(w)_{b,a}
\]
is the number of continuous lifts of any standard geometric representative of $w$
beginning at $a$ and ending at $b$; this number is independent of the representative choices.
\end{theorem}

\begin{proof}
The covering homotopy property canonically identifies the lifts obtained from
endpoint-fixed homotopic representatives of each regular factor.  Thus the count is
independent of the representative choices.  For a word containing no singular generators, ordinary path lifting is unique, so
the associated permutation matrix has entry $1$ exactly for the unique
transported endpoint and $0$ otherwise.

Consider one singular generator $\delta_i$.  Once the adjacent regular sheet states
$y_-,y_+\in A_i^{\mathrm{reg}}$ are fixed, the incidence lift theorem identifies the
possible exceptional values with
\[
A_{i,\kappa_i(y_-)}\cap A_{i,\kappa_i(y_+)}.
\]
Its cardinality is the corresponding entry of $\mathcal M_i$.

For several singular generators, choose the intermediate regular sheet
states
\[
s_1,\ldots,s_{m-1}.
\]
Ordinary regular transport between successive crossings is deterministic.
At distinct singular times in a standard representative, the exceptional values can be chosen independently once the adjacent regular states are fixed.

Hence the number of lifts realizing a specified sequence of intermediate
regular states is the product of the corresponding local incidence
multiplicities.

Summing over all possible intermediate regular states gives precisely
the matrix product
\[
P_{\eta_m}
\mathcal M_{i_m}
P_{\eta_{m-1}}
\cdots
\mathcal M_{i_1}
P_{\eta_0}.
\]
Thus its $(b,a)$ entry counts exactly the completed lifts from $a$ to
$b$.
\end{proof}

The theorem applies to isolated standard crossings; positive-duration motion inside an exceptional fibre requires additional path data internal to $F_i$.

\subsection{Multiplicity functor}
\label{subsec:multiplicity-functor}

\begin{corollary}
\label{cor:multiplicity-functor}
For finite discrete exceptional fibres,
\[
x\longmapsto A_x,
\qquad
w\longmapsto\mathscr T_{\mathbb N}(w)
\]
defines a functor
\[
\mathscr T_{\mathbb N}:
\mathcal P_{\mathrm{fc}}(X,\Sigma;\mathbf u)
\longrightarrow
\mathbf{Mat}_{\mathbb N}.
\]
\end{corollary}

\begin{proof}
The matrix assigned to a concatenation is, by definition and
Theorem~\ref{thm:multiplicity-transport}, the product of the matrices
assigned to the constituent words.  Identity words give identity
matrices.
\end{proof}

There is a support map
\[
\operatorname{supp}:
\mathbb N\longrightarrow\mathbb B,
\qquad
n\longmapsto
\begin{cases}
0,&n=0,\\
1,&n>0,
\end{cases}
\]
from the natural-number semiring to the Boolean semiring.

Applying it entrywise yields
\[
\operatorname{supp}
\bigl(
\mathscr T_{\mathbb N}(w)
\bigr)
=
\mathscr T(w).
\]

\subsection{Dependence on higher incidence}
\label{subsec:matrix-higher-incidence}

The matrix functor uses only two-sided isolated path crossings.  It therefore
depends only on the pairwise germ intersections
\[
A_{i,c}\cap A_{i,d},
\qquad c,d\in C_i.
\]

Example~\ref{ex:same-pairwise-different-higher} gives incidence data with identical
pairwise multiplicity matrices and different triple intersections.  Consequently
$\mathscr T_{\mathbb N}$ need not determine the full incidence completion.

Higher incidence is recovered by multi-arm probes in which several regular arms meet
at one defect point.

For a multi-arm probe with accumulated-end set
\[
J\subseteq C_i,
\]
the relevant multiplicity is
\[
\lambda_i(J)
=
\left|
\bigcap_{c\in J}A_{i,c}
\right|.
\]

Thus the complete finite discrete transport signature consists of the
family
\[
\lambda_i:
\mathcal P(C_i)\setminus\{\varnothing\}
\longrightarrow
\mathbb N,
\]
not merely its pairwise matrix truncation.

\section{The finite defect poset}
\label{sec:finite-defect-poset}

Sections~\ref{sec:incidence-completions}--\ref{sec:local-singular-transport}
fix the regular covering and its completion data.  In the combinatorial analysis that
follows, the regular monodromy data and the peripheral degree labels are held fixed and
suppressed: $\mathfrak D_{k,r}$ records only separation and incidence data.  Fix the
labelled exceptional set $[k]=\{1,\ldots,k\}$ and the peripheral label set
$[r]=\{1,\ldots,r\}$, with $k,r\geq1$, and vary the finite topology on $[k]$
together with its admissible endpoint incidence.

\subsection{The topology poset}
\label{subsec:topology-poset}

Let $\operatorname{Top}([k])$ be the finite set of topologies on $[k]$, ordered
in the degeneration direction by
\[
\tau\leq\sigma\quad\Longleftrightarrow\quad\sigma\subseteq\tau.
\]
Thus the discrete topology $\tau_{\mathrm{disc}}$ is minimal and the indiscrete topology
$\tau_{\mathrm{ind}}$ is maximal.

\subsection{Incidence fibres over a fixed topology}
\label{subsec:incidence-fibre-fixed-topology}

Fix $\tau\in\operatorname{Top}([k])$.

\begin{definition}[Incidence fibre]
\label{def:incidence-fibre}
Define
\[
\operatorname{Inc}_r(\tau)
\]
to be the set of all $r$-tuples
\[
\mathbf A
=
(A_1,\ldots,A_r)
\]
such that:
\begin{enumerate}[label=(\roman*)]
\item every $A_c$ is nonempty;
\item every $A_c$ is closed in $\tau$;
\item
\[
\bigcup_{c=1}^rA_c=[k].
\]
\end{enumerate}

Order $\operatorname{Inc}_r(\tau)$ coordinatewise:
\[
\mathbf A\leq\mathbf B
\quad\Longleftrightarrow\quad
A_c\subseteq B_c
\text{ for every }c.
\]
\end{definition}

Every incidence fibre has the canonical maximum
$\mathbf1_\tau:=([k],\ldots,[k])$, the full-incidence structure; minimal elements need
not be unique.

\begin{remark}
Finite unions of closed sets are closed, so $\operatorname{Inc}_r(\tau)$ is closed under coordinatewise union and is a join-semilattice.  Coordinatewise intersection can destroy nonemptiness or the covering condition, so a meet need not exist and the fibre need not be a lattice.
\end{remark}

\subsection{The total defect poset}
\label{subsec:total-defect-poset}

\begin{definition}[Finite defect poset]
\label{def:finite-defect-poset}
The \emph{finite defect poset}
\[
\mathfrak D_{k,r}
\]
is the set of all pairs
\[
(\tau,\mathbf A)
\]
with
\[
\tau\in\operatorname{Top}([k])
\]
and
\[
\mathbf A\in\operatorname{Inc}_r(\tau).
\]

Define
\[
(\tau,\mathbf A)
\leq
(\sigma,\mathbf B)
\]
if and only if
\[
\sigma\subseteq\tau
\]
and
\[
A_c\subseteq B_c
\qquad
(c=1,\ldots,r).
\]
\end{definition}

Thus a degeneration may simultaneously:

\begin{enumerate}[label=(\roman*)]
\item coarsen the exceptional topology;
\item enlarge one or more endpoint sets.
\end{enumerate}

\begin{proposition}
\label{prop:defect-poset}
The relation of
Definition~\ref{def:finite-defect-poset}
is a partial order.
\end{proposition}

\begin{proof}
Reflexivity and transitivity follow from inclusion.

If
\[
(\tau,\mathbf A)\leq(\sigma,\mathbf B)
\]
and
\[
(\sigma,\mathbf B)\leq(\tau,\mathbf A),
\]
then
\[
\sigma\subseteq\tau
\qquad\text{and}\qquad
\tau\subseteq\sigma,
\]
so
\[
\tau=\sigma.
\]
Likewise
\[
A_c=B_c
\]
for every $c$.  Hence the two elements coincide.
\end{proof}

The projection
\[
\pi_{\mathrm{sep}}:
\mathfrak D_{k,r}
\longrightarrow
\operatorname{Top}([k]),
\qquad
(\tau,\mathbf A)\longmapsto\tau
\]
will be called the \emph{separation projection}.

Its fibre over $\tau$ is exactly
\[
\operatorname{Inc}_r(\tau).
\]

\subsection{Degeneration and lift order}
\label{subsec:defect-lift-order}

Let
\[
\mathcal D=(\tau,\mathbf A)
\leq
\mathcal E=(\sigma,\mathbf B).
\]
Then $\sigma\subseteq\tau$ and $A_c\subseteq B_c$ for every $c$.

\begin{proposition}[Degeneration enlarges lift sets]
\label{prop:degeneration-enlarges-lifts}
Let
\[
\mathcal D\leq\mathcal E.
\]
For every probe with fixed regular lift, every lift admissible for
$\mathcal D$ is also admissible for $\mathcal E$.
\end{proposition}

\begin{proof}
A colouring continuous for the finer topology $\tau$ is continuous for
the coarser topology $\sigma$.

Moreover,
\[
A_c^{\mathcal D}
\subseteq
A_c^{\mathcal E}
\]
for every peripheral label.  Hence for every accumulated-end set $J$,
\[
\bigcap_{c\in J}
A_c^{\mathcal D}
\subseteq
\bigcap_{c\in J}
A_c^{\mathcal E}.
\]

Thus both continuity constraints and incidence constraints can only
weaken under degeneration.
\end{proof}

\subsection{Incidence transport under topology degeneration}
\label{subsec:incidence-transport-topology}

Let
\[
\tau\leq\sigma,
\]
so
\[
\sigma\subseteq\tau.
\]

If
\[
\mathbf A=(A_1,\ldots,A_r)
\in\operatorname{Inc}_r(\tau),
\]
the sets $A_c$ need not remain closed in the coarser topology $\sigma$.

The canonical incidence structure obtained after coarsening the topology
is therefore
\[
F_r(\tau\to\sigma)(\mathbf A)
:=
\left(
\overline{A_1}^{\,\sigma},
\ldots,
\overline{A_r}^{\,\sigma}
\right).
\]

\begin{proposition}
\label{prop:closure-incidence-map}
The tuple
\[
F_r(\tau\to\sigma)(\mathbf A)
\]
lies in
\[
\operatorname{Inc}_r(\sigma).
\]
Moreover,
\[
\mathbf A
\leq
F_r(\tau\to\sigma)(\mathbf A)
\]
coordinatewise.
\end{proposition}

\begin{proof}
Each $\sigma$-closure is nonempty and $\sigma$-closed.  Since
\[
A_c\subseteq
\overline{A_c}^{\,\sigma},
\]
their union still covers $[k]$.

The coordinatewise inequality follows from extensivity of closure.
\end{proof}

\begin{proposition}[Functoriality of closure transport]
\label{prop:closure-transport-functorial}
If
\[
\rho\subseteq\sigma\subseteq\tau,
\]
then
\[
F_r(\sigma\to\rho)
\circ
F_r(\tau\to\sigma)
=
F_r(\tau\to\rho).
\]
\end{proposition}

\begin{proof}
For every subset
\[
A\subseteq[k],
\]
one has
\[
\overline{
\overline A^{\,\sigma}
}^{\,\rho}
=
\overline A^{\,\rho},
\]
because $\rho$ is coarser than $\sigma$.

Applying this identity coordinatewise gives the result.
\end{proof}

Hence
\[
\tau
\longmapsto
\operatorname{Inc}_r(\tau)
\]
with the closure maps defines a functor from the topology degeneration
poset to finite posets.

\subsection{The one-end defect poset}
\label{subsec:r1-poset-preview}

When
\[
r=1,
\]
the incidence fibre over every topology contains only one admissible
element.

Indeed
\[
A_1\neq\varnothing
\]
and
\[
A_1=[k]
\]
are forced by endpoint-covering.

Therefore
\[
\mathfrak D_{k,1}
\cong
\operatorname{Top}([k]).
\]

After deleting the discrete and indiscrete extrema, this is the classical proper finite-topology/preorder sector considered in Section~\ref{sec:lie-boundary}.

\subsection{Symmetric-group actions}
\label{subsec:defect-poset-symmetry}

There are two natural label symmetries.

The symmetric group
\[
\mathfrak S_k
\]
acts on $[k]$ and therefore transports topologies and endpoint subsets.

The symmetric group
\[
\mathfrak S_r
\]
permutes the peripheral labels.

For
\[
(\alpha,\beta)
\in
\mathfrak S_k\times\mathfrak S_r,
\]
define
\[
(\alpha,\beta)\cdot
(\tau,A_1,\ldots,A_r)
=
\left(
\alpha_*\tau,\,
\alpha(A_{\beta^{-1}(1)}),
\ldots,
\alpha(A_{\beta^{-1}(r)})
\right),
\]
where
\[
\alpha_*\tau
=
\{\alpha(U):U\in\tau\}.
\]

\begin{proposition}
\label{prop:Skr-action}
This defines an action of
\[
\mathfrak S_k\times\mathfrak S_r
\]
on
\[
\mathfrak D_{k,r}
\]
by poset automorphisms.
\end{proposition}

\begin{proof}
Relabelling preserves nonemptiness, closedness, endpoint-covering, topology
inclusion, and coordinatewise inclusion of endpoint sets.  The group
law follows directly from composition of the two permutations.
\end{proof}

\subsection{The incidence diagram}
\label{subsec:incidence-diagram}

Let
\[
B_k:=\operatorname{Top}([k])
\]
with the degeneration order.  Proposition~\ref{prop:closure-transport-functorial}
and monotonicity of closure define a functor
\[
F_r:B_k\longrightarrow\mathbf{Poset},
\qquad
F_r(\tau)=\operatorname{Inc}_r(\tau),
\]
whose structure maps are the coordinatewise closure maps of
Section~\ref{subsec:incidence-transport-topology}.

\subsection{The Grothendieck construction}
\label{subsec:grothendieck-construction}

For a poset-valued functor
\[
F:B\to\mathbf{Poset},
\]
its Grothendieck construction
\[
\int_BF
\]
is the poset whose elements are pairs
\[
(b,x),
\qquad
x\in F(b),
\]
with order
\[
(b,x)\leq(b',x')
\]
if and only if
\[
b\leq b'
\]
and
\[
F(b\leq b')(x)\leq x'.
\]

\begin{theorem}[Grothendieck description]
\label{thm:defect-grothendieck}
There is a canonical isomorphism of posets
\[
\mathfrak D_{k,r}
\cong
\int_{B_k}F_r.
\]
\end{theorem}

\begin{proof}
The two posets have the same underlying elements:
\[
(\tau,\mathbf A),
\qquad
\mathbf A\in\operatorname{Inc}_r(\tau).
\]

In the Grothendieck construction,
\[
(\tau,\mathbf A)
\leq
(\sigma,\mathbf B)
\]
if and only if
\[
\sigma\subseteq\tau
\]
and
\[
\overline{A_c}^{\,\sigma}
\subseteq
B_c
\qquad
(c=1,\ldots,r).
\]

Since each $B_c$ is $\sigma$-closed,
\[
\overline{A_c}^{\,\sigma}
\subseteq B_c
\quad\Longleftrightarrow\quad
A_c\subseteq B_c.
\]

Hence the Grothendieck order is exactly
\[
\sigma\subseteq\tau,
\qquad
A_c\subseteq B_c
\quad(c=1,\ldots,r),
\]
which is the defining order of
\[
\mathfrak D_{k,r}.
\]
\end{proof}

\subsection{The full-incidence section}
\label{subsec:full-incidence-section}

Every incidence fibre
\[
\operatorname{Inc}_r(\tau)
\]
has the maximum
\[
\mathbf1_\tau
=
([k],\ldots,[k]).
\]

Moreover, these maxima are preserved by topology degeneration:
\[
F_r(\tau\leq\sigma)(\mathbf1_\tau)
=
\mathbf1_\sigma.
\]

Thus they define a canonical section
\[
s:
B_k\longrightarrow\mathfrak D_{k,r}
\]
given by
\[
s(\tau)
=
(\tau,\mathbf1_\tau).
\]

\begin{proposition}
\label{prop:full-incidence-section}
The map
\[
s:B_k\to\mathfrak D_{k,r}
\]
is an order embedding and satisfies
\[
\pi_{\mathrm{sep}}\circ s
=
\operatorname{id}_{B_k}.
\]
Its image is the full-incidence subposet
\[
\mathfrak F_{k,r}
:=
\left\{
(\tau,[k],\ldots,[k]):
\tau\in B_k
\right\}.
\]
\end{proposition}

\begin{proof}
If
\[
\tau\leq\sigma,
\]
then
\[
\sigma\subseteq\tau
\]
and
\[
[k]\subseteq[k]
\]
in every incidence coordinate, so
\[
s(\tau)\leq s(\sigma).
\]
Reflection of order follows immediately from the topology coordinate.
The remaining assertions follow from the definition of $s$.
\end{proof}

Hence
\[
\mathfrak F_{k,r}
\cong
B_k.
\]

\subsection{The full-incidence closure operator}
\label{subsec:full-incidence-closure}

Define
\[
\Phi:
\mathfrak D_{k,r}
\longrightarrow
\mathfrak D_{k,r}
\]
by
\[
\Phi(\tau,\mathbf A)
=
(\tau,\mathbf1_\tau).
\]

\begin{proposition}
\label{prop:Phi-closure}
The map $\Phi$ is an ascending closure operator:
\[
x\leq\Phi(x),
\qquad
\Phi^2=\Phi.
\]
Its image is exactly
\[
\mathfrak F_{k,r}=s(B_k).
\]
\end{proposition}

\begin{proof}
Extensivity follows from $A_c\subseteq[k]$ for every coordinate, and idempotence is immediate.  If
$(\tau,\mathbf A)\leq(\sigma,\mathbf B)$, then $\sigma\subseteq\tau$, hence
$(\tau,\mathbf1_\tau)\leq(\sigma,\mathbf1_\sigma)$; thus $\Phi$ is order preserving.  Its fixed points are exactly the full-incidence structures.
\end{proof}

\begin{theorem}[Full-incidence retraction]
\label{thm:full-incidence-retraction}
There is a homotopy equivalence
\[
\Delta(\mathfrak D_{k,r})
\simeq
\Delta(B_k).
\]
More precisely,
\[
\Delta(\mathfrak D_{k,r})
\simeq
\Delta(\mathfrak F_{k,r})
\cong
\Delta(\operatorname{Top}([k])).
\]
\end{theorem}

\begin{proof}
The closure operator $\Phi$ has image
\[
\mathfrak F_{k,r}.
\]
The standard closure-operator lemma for finite posets
\citep{Barmak2011} gives
\[
\Delta(\mathfrak D_{k,r})
\simeq
\Delta(\mathfrak F_{k,r}).
\]
The section $s$ identifies
\[
\mathfrak F_{k,r}
\]
with
\[
B_k.
\]
\end{proof}

Since both posets retain terminal objects, both total order complexes are contractible; the relative form below is the nontrivial statement used later.

\subsection{Restriction to a separation sector}
\label{subsec:restricted-retraction}

Let
\[
S\subseteq B_k
\]
be any subposet.

Define the corresponding restricted defect poset by
\[
\mathfrak D_{k,r}|_S
:=
\left\{
(\tau,\mathbf A)\in\mathfrak D_{k,r}:
\tau\in S
\right\}.
\]

The closure operator $\Phi$ preserves this subposet because it does not
change the topology coordinate.

\begin{theorem}[Relative full-incidence retraction]
\label{thm:relative-full-incidence-retraction}
For every subposet
\[
S\subseteq\operatorname{Top}([k]),
\]
there is a homotopy equivalence
\[
\Delta\!\left(
\mathfrak D_{k,r}|_S
\right)
\simeq
\Delta(S).
\]
\end{theorem}

\begin{proof}
Restrict $\Phi$ to
\[
\mathfrak D_{k,r}|_S.
\]
Its image is
\[
s(S)
=
\{(\tau,\mathbf1_\tau):\tau\in S\},
\]
which is isomorphic to $S$.  The same closure-operator lemma gives the result.
\end{proof}

\subsection{Incidence-blindness of total order homotopy}
\label{subsec:incidence-blindness}

Theorem~\ref{thm:relative-full-incidence-retraction} shows that, whenever all incidence structures over a topology sector $S$ are retained,
\[
\Delta(\mathfrak D_{k,r}|_S)\simeq\Delta(S).
\]
The ordinary order-complex homotopy type of the full defect poset is therefore
incidence-blind.  The incidence-sensitive constructions below instead remove the
full-incidence maximum fibrewise.

\section{The one-end sector and the partition-lattice/free-Lie comparison}
\label{sec:lie-boundary}

For $r=1$, endpoint-covering forces $A_1=[k]$, so
$\mathfrak D_{k,1}\cong\operatorname{Top}([k])$.  Delete the discrete and indiscrete
extrema and write
$\mathfrak D_{k,1}^{\circ}:=\operatorname{Top}([k])\setminus
\{\tau_{\mathrm{disc}},\tau_{\mathrm{ind}}\}$.  The homotopy type and the
representation-theoretic comparison in this section are classical.  Throughout
Sections~\ref{sec:lie-boundary} and~\ref{sec:vertical-incidence-complexes}, reduced
simplicial homology is taken in the augmented sense; in particular,
$\widetilde H_{-1}(\varnothing;R)\cong R$ and the remaining reduced homology groups
of the empty complex vanish.

\subsection{Finite topologies and preorders}
\label{subsec:topologies-preorders}

Recall that every topology $\tau$ on the finite labelled set $[k]$
determines its specialization preorder
\[
x\preceq_\tau y
\quad\Longleftrightarrow\quad
x\in\overline{\{y\}}.
\]
Conversely, every preorder on $[k]$ determines a finite topology whose
open sets are the upper sets.

Let
\[
\operatorname{Pre}([k])
\]
denote the set of all preorders on $[k]$, ordered by inclusion of
relations.

Under the degeneration order on finite topologies,
\[
\tau\leq\sigma
\quad\Longleftrightarrow\quad
\sigma\subseteq\tau,
\]
coarsening the topology enlarges the specialization preorder.

\begin{proposition}[Topology--preorder correspondence]
\label{prop:topology-preorder-poset}
The specialization construction gives an
$\mathfrak S_k$-equivariant isomorphism of posets
\[
\operatorname{Top}([k])
\cong
\operatorname{Pre}([k]),
\]
where the right-hand side is ordered by inclusion of preorder relations.

Under this correspondence:
\[
\tau_{\mathrm{disc}}
\longleftrightarrow
\Delta_k,
\]
where $\Delta_k$ is the equality preorder, and
\[
\tau_{\mathrm{ind}}
\longleftrightarrow
[k]\times[k],
\]
the universal preorder.
\end{proposition}

\begin{proof}
The finite topology--preorder correspondence is standard \citep{Alexandroff1937,Stong1966,Barmak2011}.  Under the degeneration convention used here, the order is checked as follows.

If
\[
\sigma\subseteq\tau,
\]
then every $\sigma$-open neighbourhood of a point is also constrained by
fewer opens than in $\tau$.  Consequently the specialization relation
for $\sigma$ contains that for $\tau$:
\[
\preceq_\tau
\subseteq
\preceq_\sigma.
\]
Thus degeneration order corresponds exactly to inclusion of preorder
relations.

Relabelling the underlying set commutes with both the topology and
specialization constructions, giving $\mathfrak S_k$-equivariance.
\end{proof}

Define the proper preorder poset
\[
\operatorname{Pre}_k^\circ
:=
\operatorname{Pre}([k])
\setminus
\{\Delta_k,[k]\times[k]\}.
\]

Then
\[
\mathfrak D_{k,1}^{\circ}
\cong
\operatorname{Pre}_k^\circ
\]
equivariantly.

\subsection{Directed graphs and transitive closure}
\label{subsec:preorders-directed-graphs}

Let
\[
\Omega_k
:=
[k]\times[k]
\setminus
\{(i,i):1\leq i\leq k\}.
\]

A loopless directed graph on $[k]$ is identified with a subset
\[
E\subseteq\Omega_k.
\]

Let
\[
\Delta_k^{\mathrm{NSC}}
\]
be the simplicial complex whose simplices are the edge sets of directed
graphs that are not strongly connected.

This is the complex studied by Bj\"orner and Welker
\citep{BjornerWelker1999}.

For a directed graph $D$, write
\[
\widehat D
\]
for its reflexive transitive closure.  Thus
\[
i\preceq_{\widehat D}j
\]
if and only if either $i=j$ or there is a directed path from $i$ to $j$
in $D$.

The relation $\widehat D$ is a preorder.

Moreover,
\[
D\text{ is strongly connected}
\quad\Longleftrightarrow\quad
\widehat D=[k]\times[k].
\]

If $D$ has at least one edge, then
\[
\widehat D\neq\Delta_k.
\]

Consequently transitive closure sends every nonempty face of
$\Delta_k^{\mathrm{NSC}}$ to a proper nontrivial preorder.

Let
\[
L(\Delta_k^{\mathrm{NSC}})
\]
denote the face poset of nonempty simplices of
$\Delta_k^{\mathrm{NSC}}$, ordered by inclusion.

\begin{proposition}[Transitive-closure retraction]
\label{prop:NSC-preorder-closure}
The map
\[
\operatorname{cl}:
L(\Delta_k^{\mathrm{NSC}})
\longrightarrow
\operatorname{Pre}_k^\circ,
\qquad
D\longmapsto\widehat D,
\]
is an $\mathfrak S_k$-equivariant closure map onto the proper preorder
poset.

More precisely, identifying a preorder with its set of nontrivial
ordered pairs, one has
\[
D\subseteq\widehat D,
\qquad
\widehat{\widehat D}=\widehat D.
\]
\end{proposition}

\begin{proof}
Transitive closure is extensive, idempotent, and order preserving.

If $D$ is not strongly connected, then its transitive closure is not the
universal preorder.  Since $D$ is a nonempty face, the closure is not
the equality preorder.  Hence
\[
\widehat D\in\operatorname{Pre}_k^\circ.
\]

Conversely, let
\[
P\in\operatorname{Pre}_k^\circ.
\]
Its strict relation
\[
E(P)
:=
\{(i,j):i\neq j,\ i\preceq_P j\}
\]
is a nonempty directed graph.  Since $P$ is not universal, $E(P)$ is not
strongly connected, and its reflexive transitive closure is $P$.

Thus every proper preorder lies in the image.

Finally, relabelling vertices commutes with directed paths and therefore
with transitive closure, proving equivariance.
\end{proof}

\begin{corollary}[Equivariant classical bridge]
\label{cor:equivariant-classical-bridge}
There is an $\mathfrak S_k$-equivariant homotopy equivalence
\[
\Delta(
\mathfrak D_{k,1}^{\circ}
)
\simeq_{\mathfrak S_k}
\Delta_k^{\mathrm{NSC}}.
\]
\end{corollary}

\begin{proof}
The order complex of
\[
L(\Delta_k^{\mathrm{NSC}})
\]
is the barycentric subdivision of
\[
\Delta_k^{\mathrm{NSC}},
\]
and hence is equivariantly homeomorphic to it.

Proposition~\ref{prop:NSC-preorder-closure} gives an equivariant closure
retraction from the face poset onto
\[
\operatorname{Pre}_k^\circ
\cong
\mathfrak D_{k,1}^{\circ}.
\]
The equivariant closure principle is standard; see \citet{ThevenazWebb1991}.  It gives the result.
\end{proof}

\subsection{The Bj\"orner--Welker homotopy theorem}
\label{subsec:BW-homotopy}

Björner and Welker prove that the complex of not-strongly-connected
directed graphs has the homotopy type
\[
\Delta_k^{\mathrm{NSC}}
\simeq
\bigvee^{(k-1)!}S^{2k-4}
\]
\citep{BjornerWelker1999}.

\begin{theorem}[One-end homotopy type]
\label{thm:one-end-homotopy}
For
\[
k\geq2,
\]
\[
\Delta(
\mathfrak D_{k,1}^{\circ}
)
\simeq
\bigvee^{(k-1)!}S^{2k-4}.
\]

Consequently,
\[
\widetilde H_j
\left(
\Delta(
\mathfrak D_{k,1}^{\circ}
);
\mathbb Z
\right)
=
0
\qquad
(j\neq2k-4),
\]
while
\[
\widetilde H_{2k-4}
\left(
\Delta(
\mathfrak D_{k,1}^{\circ}
);
\mathbb Z
\right)
\cong
\mathbb Z^{(k-1)!}.
\]
\end{theorem}

\subsection{The partition-lattice/free-Lie representation}
\label{subsec:free-Lie-representation}

Let $\operatorname{Lie}(k)$ denote the multilinear component of the free Lie
algebra on $k$ generators with its natural $\mathfrak S_k$-action, and let
$\overline{\Pi}_k$ denote the partition lattice on $[k]$ with its minimum and
maximum elements removed.  Under the standard convention of \citet{Stanley1982}
and \citet{Reutenauer1993},
\[
\widetilde H_{k-3}(\overline{\Pi}_k;\mathbb C)
\cong
\operatorname{Lie}(k)\otimes\operatorname{sgn}_k.
\]
The equivariant analysis of \citet{BjornerWelker1999} relates the top homology
representation of $\Delta_k^{\mathrm{NSC}}$ to this partition-lattice
representation.  No additional sign convention for the directed-graph reduction is
fixed here, and no subsequent argument depends on one.

\subsection{Persistence under incidence enrichment}
\label{subsec:Lie-persistence}

Recall the fixed separation subposet
\[
\mathfrak D_{k,r}^{\mathrm{sep}\circ}
=
\left\{
(\tau,\mathbf A):
\tau\neq\tau_{\mathrm{disc}},
\;
\tau\neq\tau_{\mathrm{ind}}
\right\}.
\]

The relative full-incidence retraction gives
\[
\Delta(
\mathfrak D_{k,r}^{\mathrm{sep}\circ}
)
\simeq
\Delta(
\mathfrak D_{k,1}^{\circ}
)
\]
for every
\[
r\geq1.
\]

The closure operator is $\mathfrak S_k$-equivariant, and the
$\mathfrak S_k$ action does not permute the peripheral labels.

\begin{corollary}[Persistence of the classical top-homology representation]
\label{cor:Lie-persists}
For every
\[
r\geq1,
\]
the $\mathfrak S_k$-module on the top reduced homology of the fixed
separation subposet satisfies
\[
\widetilde H_{2k-4}
\left(
\Delta(
\mathfrak D_{k,r}^{\mathrm{sep}\circ}
);
\mathbb C
\right)
\cong_{\mathfrak S_k}
\widetilde H_{2k-4}
\left(
\Delta_k^{\mathrm{NSC}};
\mathbb C
\right).
\]
The right-hand side is the classical representation described through the
partition-lattice/free-Lie model above.
\end{corollary}

\section{Vertical incidence complexes}
\label{sec:vertical-incidence-complexes}

Theorem~\ref{thm:relative-full-incidence-retraction} shows that the ordinary total
order-complex homotopy type does not detect the incidence fibres.

\begin{definition}[Vertical incidence complex]
\label{def:vertical-incidence-complex}
For $\tau\in\operatorname{Top}([k])$, let
\[
\mathbf1_\tau=([k],\ldots,[k])
\]
and define
\[
\mathcal V_{k,r}(\tau)
:=
\Delta
\left(
\operatorname{Inc}_r(\tau)
\setminus\{\mathbf1_\tau\}
\right).
\]
\end{definition}

When the punctured incidence fibre is empty, its order complex is taken to be
the empty simplicial complex.  The family
\[
\tau\longmapsto\mathcal V_{k,r}(\tau)
\]
is the \emph{vertical incidence profile} of the finite defect poset.
\subsection{The one-end and indiscrete extremes}
\label{subsec:vertical-extremes}

\begin{proposition}
\label{prop:vertical-r1-empty}
For every topology $\tau$,
\[
\mathcal V_{k,1}(\tau)=\varnothing.
\]
\end{proposition}

\begin{proof}
For $r=1$, endpoint-covering forces
\[
A_1=[k].
\]
Hence
\[
\operatorname{Inc}_1(\tau)
=
\{\mathbf1_\tau\}.
\]
Removing the maximum leaves the empty poset.
\end{proof}

\begin{proposition}
\label{prop:vertical-indiscrete-empty}
Let
\[
\tau_{\mathrm{ind}}
\]
be the indiscrete topology on $[k]$.  Then for every $r\geq1$,
\[
\mathcal V_{k,r}(\tau_{\mathrm{ind}})
=
\varnothing.
\]
\end{proposition}

\begin{proof}
The only nonempty closed subset of an indiscrete space is $[k]$ itself.
Hence every endpoint set is forced to equal $[k]$, so
\[
\operatorname{Inc}_r(\tau_{\mathrm{ind}})
=
\{\mathbf1_{\tau_{\mathrm{ind}}}\}.
\]
\end{proof}

\subsection{Order-theoretic description for a finite topology}
\label{subsec:vertical-order-description}

Let
\[
\rho:[k]\longrightarrow P
\]
be the Kolmogorov quotient of the finite topology $\tau$, with
specialization poset $P$.

Every closed subset of $[k]$ is saturated for topological
indistinguishability and is therefore the inverse image of a unique
lower set in $P$.

Let
\[
\mathcal L(P)
\]
denote the lattice of lower sets of $P$, and write
\[
\mathcal L^\times(P)
:=
\mathcal L(P)\setminus\{\varnothing\}.
\]

\begin{definition}[Order-incidence fibre]
\label{def:order-incidence-fibre}
Define
\[
\operatorname{Inc}_r(P)
\]
to be the poset of tuples
\[
(L_1,\ldots,L_r)
\in
\bigl(\mathcal L^\times(P)\bigr)^r
\]
satisfying
\[
\bigcup_{c=1}^rL_c=P,
\]
ordered coordinatewise by inclusion.
\end{definition}

\begin{proposition}
\label{prop:vertical-depends-T0}
There is a canonical isomorphism of posets
\[
\operatorname{Inc}_r(\tau)
\cong
\operatorname{Inc}_r(P).
\]
Consequently,
\[
\mathcal V_{k,r}(\tau)
\cong
\Delta
\left(
\operatorname{Inc}_r(P)
\setminus\{(P,\ldots,P)\}
\right).
\]
In particular, the vertical incidence complex depends only on the
Kolmogorov quotient poset $P$ and not on the multiplicities of its
indistinguishability classes.
\end{proposition}

\begin{proof}
A nonempty $\tau$-closed subset $A\subseteq[k]$ has the form
\[
A=\rho^{-1}(L)
\]
for a unique nonempty lower set
\[
L\subseteq P.
\]
Moreover,
\[
\bigcup_{c=1}^rA_c=[k]
\]
if and only if
\[
\bigcup_{c=1}^rL_c=P.
\]
The correspondence preserves coordinatewise inclusion and sends the
full-incidence maximum to
\[
(P,\ldots,P).
\]
\end{proof}

\subsection{The discrete topology and edge-cover posets}
\label{subsec:discrete-edge-covers}

Now let
\[
\tau=\tau_{\mathrm{disc}}.
\]
Every subset of $[k]$ is closed.

An incidence tuple
\[
(A_1,\ldots,A_r)
\]
can be identified with a bipartite graph
\[
G\subseteq K_{r,k}
\]
whose left vertex class is $[r]$, whose right vertex class is $[k]$, and
whose edge set is
\[
E(G)
=
\{(c,a):a\in A_c\}.
\]

The conditions
\[
A_c\neq\varnothing
\]
and
\[
\bigcup_cA_c=[k]
\]
say exactly that $G$ has no isolated vertices.

\begin{proposition}
\label{prop:discrete-inc-edge-cover}
The poset
\[
\operatorname{Inc}_r(\tau_{\mathrm{disc}})
\]
is canonically isomorphic to the poset of spanning edge covers of
\[
K_{r,k},
\]
ordered by inclusion of edge sets.
\end{proposition}

The full-incidence maximum corresponds to the complete bipartite graph
\[
K_{r,k}.
\]

Therefore
\[
\mathcal V_{k,r}(\tau_{\mathrm{disc}})
\]
is the order complex of all proper edge covers of $K_{r,k}$.

\subsection{Complementation and bounded-degree graphs}
\label{subsec:complement-bounded-degree}

Let
\[
\Omega=[r]\times[k]
\]
be the edge set of $K_{r,k}$.

For an incidence graph $G$, define its edge complement
\[
D
:=
\Omega\setminus E(G).
\]

The requirement that every left vertex $c\in[r]$ be incident with at
least one edge of $G$ becomes
\[
\deg_D(c)\leq k-1.
\]

Likewise, endpoint-covering says every
$a\in[k]$ has at least one incident edge in $G$, which becomes
\[
\deg_D(a)\leq r-1.
\]

Thus define the simplicial complex
\[
\mathcal C_{k,r}
\]
on vertex set
\[
[r]\times[k]
\]
by
\[
D\in\mathcal C_{k,r}
\quad\Longleftrightarrow\quad
\begin{cases}
\deg_D(c)\leq k-1
& (c\in[r]),\\[1mm]
\deg_D(a)\leq r-1
& (a\in[k]).
\end{cases}
\]

The degree bounds are hereditary under deletion of edges, so
\[
\mathcal C_{k,r}
\]
is indeed an abstract simplicial complex.

Equivalently, $\mathcal C_{k,r}$ is a bipartite bounded-degree graph complex in the classical framework of \citet{KaraguezianReinerWachs2001}; in the chessboard language, compare also \citet{JojicVrecicaZivaljevic2017}.

\begin{theorem}[Discrete vertical fibre as a generalized chessboard complex]
\label{thm:discrete-vertical-chessboard}
There is a canonical homeomorphism
\[
\mathcal V_{k,r}(\tau_{\mathrm{disc}})
\cong
\operatorname{sd}(\mathcal C_{k,r}),
\]
where
\[
\operatorname{sd}
\]
denotes barycentric subdivision.

Consequently,
\[
\mathcal V_{k,r}(\tau_{\mathrm{disc}})
\simeq
\mathcal C_{k,r}.
\]
\end{theorem}

\begin{proof}
Complementation gives an order-reversing bijection between:

\[
\operatorname{Inc}_r(\tau_{\mathrm{disc}})
\setminus
\{\mathbf1_{\tau_{\mathrm{disc}}}\}
\]
and the poset of nonempty simplices of
\[
\mathcal C_{k,r}.
\]

Indeed, removing the full-incidence graph means that its complement
$D$ is nonempty.

The order complex of a poset is canonically isomorphic to that of its
opposite poset, while the order complex of the nonempty face poset of a
simplicial complex is its barycentric subdivision.

Hence
\[
\mathcal V_{k,r}(\tau_{\mathrm{disc}})
\cong
\operatorname{sd}(\mathcal C_{k,r}).
\]
\end{proof}

For the discrete topology, the vertical fibre is therefore an established bounded-degree graph complex; the remaining question is how this fibre changes under topology coarsening.

\subsection{Collapse to full incidence under topology degeneration}
\label{subsec:vertical-annihilation}

Let
\[
\tau\leq\sigma,
\]
so
\[
\sigma\subseteq\tau.
\]

The closure transport
\[
F_r(\tau\leq\sigma):
\operatorname{Inc}_r(\tau)
\longrightarrow
\operatorname{Inc}_r(\sigma)
\]
sends
\[
\mathbf A=(A_1,\ldots,A_r)
\]
to
\[
\overline{\mathbf A}^{\,\sigma}
=
\left(
\overline{A_1}^{\,\sigma},
\ldots,
\overline{A_r}^{\,\sigma}
\right).
\]

Although maxima map to maxima, a nonmaximal incidence tuple may also
become full incidence after closure.

This prevents the punctured fibres
\[
\operatorname{Inc}_r(\tau)\setminus\{\mathbf1_\tau\}
\]
from forming an ordinary functor over all of
\[
\operatorname{Top}([k]).
\]

\begin{definition}[Full-incidence collapse locus]
\label{def:vertical-annihilation}
For
\[
\tau\leq\sigma,
\]
define
\[
\mathcal A_r(\tau,\sigma)
:=
\left\{
\mathbf A\in
\operatorname{Inc}_r(\tau)\setminus\{\mathbf1_\tau\}:
F_r(\tau\leq\sigma)(\mathbf A)
=
\mathbf1_\sigma
\right\}.
\]
\end{definition}

These are precisely the proper incidence structures that become
full-incidence after the topology is coarsened from $\tau$ to $\sigma$.

Define the surviving subposet
\[
\operatorname{Surv}_r(\tau,\sigma)
:=
\left(
\operatorname{Inc}_r(\tau)\setminus\{\mathbf1_\tau\}
\right)
\setminus
\mathcal A_r(\tau,\sigma).
\]

Closure therefore restricts to an order-preserving map
\[
F_r^{\mathrm{surv}}(\tau\leq\sigma):
\operatorname{Surv}_r(\tau,\sigma)
\longrightarrow
\operatorname{Inc}_r(\sigma)
\setminus\{\mathbf1_\sigma\}.
\]

\subsection{Characterizing collapse to full incidence}
\label{subsec:annihilation-criterion}

\begin{proposition}[Full-incidence collapse criterion]
\label{prop:annihilation-criterion}
Let
\[
\tau\leq\sigma
\]
and
\[
\mathbf A=(A_1,\ldots,A_r)
\in
\operatorname{Inc}_r(\tau)
\setminus\{\mathbf1_\tau\}.
\]

Then
\[
\mathbf A\in\mathcal A_r(\tau,\sigma)
\]
if and only if
\[
\overline{A_c}^{\,\sigma}
=
[k]
\qquad
(c=1,\ldots,r).
\]
\end{proposition}

\begin{proof}
The image under closure transport is
\[
F_r(\tau\leq\sigma)(\mathbf A)
=
\left(
\overline{A_1}^{\,\sigma},
\ldots,
\overline{A_r}^{\,\sigma}
\right).
\]
This equals
\[
\mathbf1_\sigma
=
([k],\ldots,[k])
\]
exactly when every coordinate closure is the full space.
\end{proof}

\subsection{The vertical homology profile}
\label{subsec:vertical-homology-profile}

Define
\[
\mathsf V_{k,r}(\tau)
:=
\widetilde H_*
\left(
\mathcal V_{k,r}(\tau);
\mathbb Z
\right).
\]

The collection
\[
\left\{
\mathsf V_{k,r}(\tau)
\right\}_{\tau\in\operatorname{Top}([k])}
\]
will be called the \emph{vertical homology profile}.

The profile records the reduced homology of the punctured incidence fibre
separately at each separation level.

Closure may send vertices of the punctured source fibre to the deleted maximum, so it does not induce a simplicial map on the entire punctured fibre.

\section{Separation--incidence grading and obstruction to a horizontal--vertical bicomplex}
\label{sec:separation-incidence-grading}

The Hasse diagram distinguishes topology changes from incidence changes, but this edge-colouring does not induce a canonical bicomplex on the ordinary order-complex chains.

\subsection{Horizontal and vertical Hasse edges}
\label{subsec:horizontal-vertical-Hasse}

A cover relation in $\mathfrak D_{k,r}$ is called \emph{vertical} if its topology
coordinate is fixed and \emph{horizontal} otherwise.

\begin{proposition}[Hasse-edge separation]
\label{prop:Hasse-bicoloring}
If
\[
(\tau,\mathbf A)\lessdot(\sigma,\mathbf B)
\]
in $\mathfrak D_{k,r}$, then exactly one of the following holds:
\begin{enumerate}[label=(\roman*)]
\item $\tau=\sigma$ and $\mathbf A\lessdot\mathbf B$ in
$\operatorname{Inc}_r(\tau)$;
\item $\tau\neq\sigma$ and $\mathbf A=\mathbf B$.
\end{enumerate}
Hence every Hasse edge has a unique vertical or horizontal type.
\end{proposition}

\begin{proof}
If $\tau=\sigma$, minimality is exactly minimality inside the incidence fibre.  Suppose
$\tau\neq\sigma$.  Since
$(\tau,\mathbf A)<(\sigma,\mathbf B)$, one has $\sigma\subsetneq\tau$ and
$A_c\subseteq B_c$ for all $c$.  Each $B_c$ is $\sigma$-closed; because every
$\sigma$-open set is $\tau$-open, each $B_c$ is also $\tau$-closed.  Thus
$\mathbf B\in\operatorname{Inc}_r(\tau)$.  If $\mathbf A\neq\mathbf B$, then
\[
(\tau,\mathbf A)<(\tau,\mathbf B)<(\sigma,\mathbf B),
\]
contradicting that the original relation is a cover.  Therefore
$\mathbf A=\mathbf B$.
\end{proof}

Closure enlargement under topology coarsening remains relevant for arbitrary arrows
$\tau\leq\sigma$, but it does not occur across a horizontal Hasse edge.

\subsection{Move counts on saturated chains}
\label{subsec:move-counts}

Let
\[
\mathcal C:x_0\lessdot x_1\lessdot\cdots\lessdot x_N
\]
be a saturated chain.  Let $p(\mathcal C)$ and $q(\mathcal C)$ denote, respectively,
the numbers of horizontal and vertical cover relations in $\mathcal C$; then
$p(\mathcal C)+q(\mathcal C)=N$.

The ordered pair
\[
(p(\mathcal C),q(\mathcal C))
\]
records the number of separation and incidence moves in the chosen
saturated degeneration chain.

This assigns a bidegree to each \emph{saturated chain}; it does not define a bigrading of the simplicial chain complex of the order complex.

\subsection{Why the move count is not an invariant of an interval}
\label{subsec:move-count-not-invariant}

For a graded poset, all saturated chains in a fixed interval have the
same length.  No corresponding statement holds for the full
finite-topology poset.

Consequently two saturated chains with the same endpoints may have
different lengths.  Even when their total lengths agree, there is no
general reason for the horizontal and vertical counts to coincide.

The assignment
\[
[x,y]\longmapsto(p,q)
\]
is therefore not well defined in general.

\begin{remark}
This prevents a canonical rank function of the form
\[
\operatorname{rk}(x)
=
(\operatorname{rk}_{\mathrm{sep}}(x),
\operatorname{rk}_{\mathrm{inc}}(x))
\]
on the complete defect poset.  Neither component has a canonical
global rank without imposing additional restrictions.
\end{remark}

Accordingly, the pair
\[
(p,q)
\]
belongs to a chosen saturated path in the Hasse diagram, not to its
endpoints.

\begin{example}[Two saturated chains with different move counts]
\label{ex:different-move-counts}
Let $k=3$, $r=2$, and write $[3]=\{1,2,3\}$.  Consider
\[
\tau=\{\varnothing,\{1\},\{1,2\},[3]\},
\qquad
\rho=\{\varnothing,[3]\},
\]
and the intermediate topologies
\[
\sigma=\{\varnothing,\{1,2\},[3]\},
\qquad
\sigma'=\{\varnothing,\{1\},[3]\}.
\]
Set
\[
\mathbf A=([3],\{3\}),
\qquad
\mathbf A'=([3],\{2,3\}),
\qquad
\mathbf B=([3],[3]).
\]
All displayed endpoint sets are closed in the topology with which they are paired.  The interval from $(\tau,\mathbf A)$ to $(\rho,\mathbf B)$ contains the saturated chain
\[
(\tau,\mathbf A)
\lessdot
(\sigma,\mathbf A)
\lessdot
(\sigma,\mathbf B)
\lessdot
(\rho,\mathbf B),
\]
whose edge types are horizontal, vertical, horizontal and hence have move count $(2,1)$.  It also contains the saturated chain
\[
(\tau,\mathbf A)
\lessdot
(\tau,\mathbf A')
\lessdot
(\sigma',\mathbf A')
\lessdot
(\sigma',\mathbf B)
\lessdot
(\rho,\mathbf B),
\]
whose edge types are vertical, horizontal, vertical, horizontal and hence have move count $(2,2)$.  Each displayed relation is a cover: the vertical steps are covers in the corresponding incidence fibres, and the horizontal steps admit no topology strictly between the displayed topologies for which the fixed incidence tuple remains admissible.  Thus the move count is not an invariant of an interval, already in $\mathfrak D_{3,2}$.
\end{example}

\subsection{Failure of the simplicial boundary to preserve move type}
\label{subsec:boundary-obstruction}

A simplex of the order complex is an arbitrary strict chain
\[
x_0<x_1<\cdots<x_n.
\]
Its simplicial boundary is
\[
\partial[x_0<\cdots<x_n]
=
\sum_{j=0}^n
(-1)^j
[x_0<\cdots<\widehat{x_j}<\cdots<x_n].
\]

The adjacent relations
\[
x_{j-1}<x_j<x_{j+1}
\]
need not be cover relations.

Even if a simplex is represented by a saturated chain, deleting the
intermediate element $x_j$ replaces two Hasse moves by the composite
relation
\[
x_{j-1}<x_{j+1}.
\]

If the two removed edges have different colours, the composite relation
is neither a single horizontal Hasse edge nor a single vertical Hasse
edge.

Hence the ordinary simplicial boundary does not split canonically as
\[
\partial_{\mathrm{sep}}
+
\partial_{\mathrm{inc}}
\]
by deleting edges of different colours.

\begin{proposition}[Horizontal--vertical bicomplex obstruction]
\label{prop:naive-bicomplex-obstruction}
The canonical horizontal--vertical colouring of the Hasse diagram does
not, in general, induce a canonical bicomplex structure on
\[
C_*
\bigl(
\Delta(\mathfrak D_{k,r})
\bigr)
\]
by counting horizontal and vertical cover relations.
\end{proposition}

\begin{proof}
A bicomplex grading would require the simplicial boundary to decompose
into operators changing one grading at a time.

However, simplices of the order complex are chains rather than Hasse
paths, and deleting a vertex composes two adjacent order relations.
The resulting composite need not inherit either Hasse colour.

Moreover, the number and types of Hasse edges refining a non-cover
relation need not be independent of the chosen saturated refinement.

Therefore the proposed bidegree is not intrinsic to simplicial chains
and the required decomposition of the boundary is not canonically
defined.
\end{proof}

\subsection{Survival data under topology degeneration}
\label{subsec:survival-data-horizontal}

For
\[
\tau\leq\sigma,
\]
define
\[
\operatorname{Surv}_r(\tau,\sigma)
\]
to be the proper incidence structures over $\tau$ that do not map to
full incidence under closure.

The closure map induces
\[
F_r^{\mathrm{surv}}(\tau\leq\sigma):
\operatorname{Surv}_r(\tau,\sigma)
\longrightarrow
\operatorname{Inc}_r(\sigma)
\setminus
\{\mathbf1_\sigma\}.
\]

Taking order complexes yields a simplicial map
\[
\Delta\operatorname{Surv}_r(\tau,\sigma)
\longrightarrow
\mathcal V_{k,r}(\sigma).
\]

The inclusion
\[
\operatorname{Surv}_r(\tau,\sigma)
\hookrightarrow
\operatorname{Inc}_r(\tau)
\setminus
\{\mathbf1_\tau\}
\]
gives
\[
\Delta\operatorname{Surv}_r(\tau,\sigma)
\longrightarrow
\mathcal V_{k,r}(\tau).
\]

Every topology degeneration determines a span
\[
\mathcal V_{k,r}(\tau)
\longleftarrow
\Delta\operatorname{Surv}_r(\tau,\sigma)
\longrightarrow
\mathcal V_{k,r}(\sigma).
\]

The middle term records exactly those proper incidence structures that survive the
topology degeneration.

\subsection{The vertical correspondence diagram}
\label{subsec:vertical-correspondence-diagram}

For every topology arrow
\[
\tau\leq\sigma
\]
associate the span
\[
\mathcal V_{k,r}(\tau)
\longleftarrow
\mathcal S_{k,r}(\tau,\sigma)
\longrightarrow
\mathcal V_{k,r}(\sigma),
\]
where
\[
\mathcal S_{k,r}(\tau,\sigma)
:=
\Delta
\operatorname{Surv}_r(\tau,\sigma).
\]

Composition of topology degenerations
\[
\tau\leq\sigma\leq\rho
\]
requires comparing
\[
\operatorname{Surv}_r(\tau,\rho)
\]
with the incidence structures that survive first to $\sigma$ and then
to $\rho$.

Because closure transport is functorial,
\[
F_r(\tau\leq\rho)
=
F_r(\sigma\leq\rho)
\circ
F_r(\tau\leq\sigma),
\]
an incidence structure survives directly to $\rho$ if and only if its
image at $\sigma$ survives from $\sigma$ to $\rho$.

\begin{proposition}[Survival under composition]
\label{prop:survival-composition}
For
\[
\tau\leq\sigma\leq\rho,
\]
one has
\[
\operatorname{Surv}_r(\tau,\rho)
=
\left\{
\mathbf A\in
\operatorname{Surv}_r(\tau,\sigma):
F_r(\tau\leq\sigma)(\mathbf A)
\in
\operatorname{Surv}_r(\sigma,\rho)
\right\}.
\]
\end{proposition}

\begin{proof}
By definition,
\[
\mathbf A\in
\operatorname{Surv}_r(\tau,\rho)
\]
if and only if
\[
F_r(\tau\leq\rho)(\mathbf A)
\neq
\mathbf1_\rho.
\]

Functoriality gives
\[
F_r(\tau\leq\rho)(\mathbf A)
=
F_r(\sigma\leq\rho)
\left(
F_r(\tau\leq\sigma)(\mathbf A)
\right).
\]

If the first image
\[
F_r(\tau\leq\sigma)(\mathbf A)
\]
were already full incidence, its further image would remain full
incidence.  Hence direct survival implies survival to $\sigma$.

Under that condition, direct survival to $\rho$ is precisely survival of
the intermediate image from $\sigma$ to $\rho$.
\end{proof}

\section{Conclusion and open problems}
\label{sec:conclusion}

A finite covering of $X\setminus\Sigma$ determines peripheral monodromy.  An incidence
completion supplements it by exceptional fibres and a closed many-to-many relation between peripheral
component germs and exceptional points, and the tail-saturated topology realizes that relation
intrinsically.  Endpoint intersections then govern lift extension, discrete local homology,
reconstruction from probes, and singular transport.  For finite exceptional fibres the datum reduces
to $(P,\mu,I)$ and admits the structured monodromy, automorphism, and quotient descriptions of
Sections~\ref{sec:structured-monodromy} and~\ref{sec:singular-automorphisms}.

The finite defect poset organizes the variation of the exceptional topology and incidence.  Its
full-incidence closure retraction shows that ordinary total order-complex homotopy retains only the
finite-topology direction.  The one-end sector is classical.  The vertical complexes and their survival spans
retain incidence-dependent information discarded by the ordinary total order-complex homotopy type.

Open problems include the computation of $\widetilde H_*(\mathcal V_{k,r}(\tau))$ across
specialization types, a chain or correspondence theory compatible with the survival spans, transport
for positive-dimensional defect strata, and inverse-limit versions of the finite structured theory.

\bibliographystyle{plainnat}
\bibliography{references}

@book{Hatcher2002,
  author    = {Hatcher, Allen},
  title     = {Algebraic Topology},
  publisher = {Cambridge University Press},
  address   = {Cambridge},
  year      = {2002}
}

@book{Spanier1966,
  author    = {Spanier, Edwin H.},
  title     = {Algebraic Topology},
  publisher = {McGraw-Hill},
  address   = {New York},
  year      = {1966}
}

@book{Brown2006,
  author    = {Brown, Ronald},
  title     = {Topology and Groupoids},
  edition   = {3},
  publisher = {BookSurge},
  address   = {Charleston, SC},
  year      = {2006}
}

@article{Woolf2009,
  author  = {Woolf, Jonathan},
  title   = {The fundamental category of a stratified space},
  journal = {Journal of Homotopy and Related Structures},
  volume  = {4},
  number  = {1},
  pages   = {359--387},
  year    = {2009},
  eprint  = {0811.2580},
  archivePrefix = {arXiv},
  primaryClass  = {math.AT}
}

@article{CalcutGompfMcCarthy2012,
  author  = {Calcut, Jack S. and Gompf, Robert E. and McCarthy, John D.},
  title   = {On fundamental groups of quotient spaces},
  journal = {Topology and its Applications},
  volume  = {159},
  number  = {1},
  pages   = {322--330},
  year    = {2012},
  doi     = {10.1016/j.topol.2011.09.038},
  eprint  = {0904.4465},
  archivePrefix = {arXiv},
  primaryClass  = {math.GN}
}

@article{Alexandroff1937,
  author  = {Alexandroff, P.},
  title   = {Diskrete R{\"a}ume},
  journal = {Matematicheskii Sbornik},
  series  = {N.S.},
  volume  = {2(44)},
  pages   = {501--518},
  year    = {1937}
}

@article{Stong1966,
  author  = {Stong, Robert E.},
  title   = {Finite topological spaces},
  journal = {Transactions of the American Mathematical Society},
  volume  = {123},
  pages   = {325--340},
  year    = {1966},
  doi     = {10.1090/S0002-9947-1966-0195042-2}
}

@article{McCord1966,
  author  = {McCord, Michael C.},
  title   = {Singular homology groups and homotopy groups of finite topological spaces},
  journal = {Duke Mathematical Journal},
  volume  = {33},
  number  = {3},
  pages   = {465--474},
  year    = {1966},
  doi     = {10.1215/S0012-7094-66-03352-7}
}

@book{Barmak2011,
  author    = {Barmak, Jonathan A.},
  title     = {Algebraic Topology of Finite Topological Spaces and Applications},
  series    = {Lecture Notes in Mathematics},
  volume    = {2032},
  publisher = {Springer},
  address   = {Berlin},
  year      = {2011},
  doi       = {10.1007/978-3-642-22003-6}
}

@article{Niefield1982,
  author  = {Niefield, Susan B.},
  title   = {Cartesianness: Topological Spaces, Uniform Spaces, and Affine Schemes},
  journal = {Journal of Pure and Applied Algebra},
  volume  = {23},
  number  = {2},
  pages   = {147--167},
  year    = {1982},
  doi     = {10.1016/0022-4049(82)90004-4}
}

@article{Niefield2012,
  author  = {Niefield, Susan B.},
  title   = {The Glueing Construction and Double Categories},
  journal = {Journal of Pure and Applied Algebra},
  volume  = {216},
  number  = {8--9},
  pages   = {1827--1836},
  year    = {2012},
  doi     = {10.1016/j.jpaa.2012.02.021},
  eprint  = {1112.1307},
  archivePrefix = {arXiv},
  primaryClass  = {math.CT}
}

@article{deSouza2023,
  author  = {de Souza, Lucas H. R.},
  title   = {Glueing spaces without identifying points},
  journal = {Topology and its Applications},
  volume  = {338},
  pages   = {108671},
  year    = {2023},
  doi     = {10.1016/j.topol.2023.108671},
  eprint  = {2004.01845},
  archivePrefix = {arXiv},
  primaryClass  = {math.GN}
}

@incollection{Fox1957,
  author    = {Fox, Ralph H.},
  title     = {Covering spaces with singularities},
  booktitle = {Algebraic Geometry and Topology: A Symposium in Honor of S. Lefschetz},
  pages     = {243--257},
  publisher = {Princeton University Press},
  address   = {Princeton, NJ},
  year      = {1957}
}

@article{Montesinos2005,
  author  = {Montesinos-Amilibia, Jos{\'e} Mar{\'i}a},
  title   = {Branched coverings after {Fox}},
  journal = {Bolet{\'i}n de la Sociedad Matem{\'a}tica Mexicana},
  series  = {3},
  volume  = {11},
  number  = {1},
  pages   = {19--64},
  year    = {2005}
}

@article{BjornerWelker1999,
  author  = {Bj{\"o}rner, Anders and Welker, Volkmar},
  title   = {Complexes of Directed Graphs},
  journal = {SIAM Journal on Discrete Mathematics},
  volume  = {12},
  number  = {4},
  pages   = {413--424},
  year    = {1999},
  doi     = {10.1137/S0895480198338724}
}

@article{Stanley1982,
  author  = {Stanley, Richard P.},
  title   = {Some aspects of groups acting on finite posets},
  journal = {Journal of Combinatorial Theory, Series A},
  volume  = {32},
  pages   = {132--161},
  year    = {1982},
  doi     = {10.1016/0097-3165(82)90017-6}
}

@book{Reutenauer1993,
  author    = {Reutenauer, Christophe},
  title     = {Free Lie Algebras},
  series    = {London Mathematical Society Monographs, New Series},
  volume    = {7},
  publisher = {Oxford University Press},
  address   = {Oxford},
  year      = {1993}
}

@article{JojicVrecicaZivaljevic2017,
  author  = {Joji{\'c}, Du{\v{s}}ko and Vre{\'c}ica, Sini{\v{s}}a T. and {\v{Z}}ivaljevi{\'c}, Rade T.},
  title   = {Multiple chessboard complexes and the colored Tverberg problem},
  journal = {Journal of Combinatorial Theory, Series A},
  volume  = {145},
  pages   = {400--425},
  year    = {2017},
  doi     = {10.1016/j.jcta.2016.08.008}
}

@article{KaraguezianReinerWachs2001,
  author  = {Karaguezian, Dikran B. and Reiner, Victor and Wachs, Michelle L.},
  title   = {Matching Complexes, Bounded Degree Graph Complexes, and Weight Spaces of {$GL_n$}-Complexes},
  journal = {Journal of Algebra},
  volume  = {239},
  number  = {1},
  pages   = {77--92},
  year    = {2001},
  doi     = {10.1006/jabr.2000.8654}
}

@book{McCleary2001,
  author    = {McCleary, John},
  title     = {A User's Guide to Spectral Sequences},
  edition   = {2},
  series    = {Cambridge Studies in Advanced Mathematics},
  volume    = {58},
  publisher = {Cambridge University Press},
  address   = {Cambridge},
  year      = {2001},
  doi       = {10.1017/CBO9780511626289}
}

@article{ThevenazWebb1991,
  author  = {Th{\'e}venaz, Jacques and Webb, Peter J.},
  title   = {Homotopy Equivalence of Posets with a Group Action},
  journal = {Journal of Combinatorial Theory, Series A},
  volume  = {56},
  number  = {2},
  pages   = {173--181},
  year    = {1991},
  doi     = {10.1016/0097-3165(91)90030-K}
}

\end{document}